\documentclass[aos,citesort,dvips]{arximspdf}
\usepackage{mathbh}
\usepackage[mathscr]{euscript}

\doi{10.1214/10-AOS843}
\volume{39}
\issue{1}
\pubyear{2011}
\firstpage{333}
\lastpage{361}

\makeatletter

\newcommand{\er}{\mathit{er}}

\newcommand{\UB}{\mathit{UB}}
\newcommand{\LB}{\mathit{LB}}

\newcommand{\D}{\mathcal D}
\newcommand{\X}{\mathcal X}
\newcommand{\Y}{\mathcal Y}
\newcommand{\LL}{\mathcal L}
\newcommand{\C}{\mathbb C}
\newcommand{\PP}{\mathbb P}
\newcommand{\nats}{\mathbb{N}}
\newcommand{\E}{\mathbb E}
\newcommand{\Data}{\mathcal Z}

\newcommand{\bound}{\mathscr{E}}

\newcommand{\diam}{\operatorname{diam}}

\newsavebox{\savepar}
\newenvironment{bigboxit}{\begin{center}\begin{lrbox}{\savepar}
\begin{minipage}[h]{4.8in}
\normalfont
\begin{flushleft}}
{\end{flushleft}\end{minipage}\end{lrbox}\fbox{\usebox{\savepar}}
\end{center}}

\newtheorem{theorem}{Theorem}
\newtheorem{lemma}{Lemma}
\newproclaim{definition}{Definition}
\newproclaim{condition}{Condition}

\makeatother

\begin{document}
\begin{frontmatter}

\title{Rates of convergence in active learning}
\runtitle{Rates of convergence in active learning}

\begin{aug}
\author[A]{\fnms{Steve} \snm{Hanneke}\corref{}\thanksref{t1}\ead[label=e1]{shanneke@stat.cmu.edu}}
\runauthor{S. Hanneke}
\affiliation{Carnegie Mellon University}
\address[A]{Department of Statistics \\
Carnegie Mellon University \\
5000 Forbes Avenue \\
Pittsburgh, Pennsylvania 15213\\
USA \\
\printead{e1}} %adresu isvedimo komanda gale!
\end{aug}

\thankstext{t1}{Supported by the NSF Grant IIS-0713379 and an IBM Ph.D.
Fellowship.}

% HISTORY:
\received{\smonth{8} \syear{2009}}
\revised{\smonth{6} \syear{2010}}

% ABSTRACT
%
\begin{abstract}
We study the rates of convergence in generalization error achievable by
active learning under various types of label noise.
Additionally, we study the general problem of model selection for
active learning with a nested hierarchy of hypothesis classes and
propose an algorithm whose error rate provably converges to the best
achievable error among classifiers in the hierarchy at a rate adaptive
to both the complexity of the optimal classifier and the noise
conditions. In particular, we state sufficient conditions for these
rates to be dramatically faster than those achievable by passive
learning.
\end{abstract}

% KEYWORDS
%

\end{frontmatter}

%s1 ###
\section{Introduction}

\textit{Active learning} refers to a family of powerful supervised
learning protocols capable of producing more accurate classifiers while
using a smaller number of labeled data points than traditional (passive)
learning methods.
Here we study a variant known as \textit{pool-based} active learning,
in which a learning algorithm is given
access to a large pool of unlabeled data (i.e., only the covariates are
visible),
and is allowed to sequentially request the label (response variable) of
any particular data points from that pool. The
objective is to learn a function that accurately predicts the labels of new
points, while minimizing the number of label requests.
Thus, this is a type of sequential design scenario for a function
estimation problem.
This contrasts with passive learning, where the labeled data are
sampled at random.
In comparison, by more carefully selecting which points should be labeled,
%focusing on the more informative examples for the learning problem,
active learning can often significantly decrease the total amount of effort
required for data annotation.
This can be particularly interesting for tasks where unlabeled
data are available in abundance, but label information comes
only through significant effort or cost.

Recently, there have been a series of exciting advances on the topic
of active learning with arbitrary classification noise (the so-called
\textit{agnostic} PAC model \cite{KSS94}), resulting in several new algorithms
capable of achieving improved convergence rates compared to passive
learning under certain conditions. The first, proposed
by Balcan, Beygelzimer and Langford \cite{balcan06} was the $A^2$
(agnostic active) algorithm, which
provably never has significantly worse rates of convergence than
passive learning by empirical risk minimization. This algorithm was
later analyzed in detail in \cite{hanneke07b}, where it was
found that a complexity measure called the \textit{disagreement
coefficient} characterizes the worst-case convergence rates achieved
by $A^2$ for any given hypothesis class, data distribution and best
achievable error rate in the class. The next major advance was
by Dasgupta, Hsu and Monteleoni \cite{dasgupta07}, who proposed a new
algorithm, and proved that
it improves the dependence of the convergence rates on the
disagreement coefficient compared to $A^2$. Both algorithms are
defined below in Section \ref{sec:algorithms}. While all of these
advances are encouraging, they are limited in two ways. First, the
convergence rates that have been proven for these algorithms typically
only improve the dependence on the magnitude of the noise (more
precisely, the noise rate of the hypothesis class), compared to
passive learning. Thus, in an asymptotic sense, for nonzero noise
rates these results represent at best a constant factor improvement
over passive learning. Second, these results are limited to learning
with a fixed hypothesis class of limited expressiveness, so that
convergence to the Bayes error rate is not always a possibility.

%%% I don't mention their work on boundary fragment classes,
%%% since the noise conditions there are much more restrictive than
%Tsybakov noise,
%%% and I don't analyze those more restrictive assumptions in this work
%at all
On the first of these limitations, recent work by Castro and Nowak
\cite{castro08} on learning threshold classifiers discovered
that if certain parameters of the noise distribution are \textit{known}
(namely, parameters related to Tsybakov's margin
conditions), then we can achieve strict improvements in the
asymptotic convergence rate
%dependence on those noise parameters
via a specific active learning algorithm designed to take advantage of
that knowledge for thresholds.
Subsequently, Balcan, Broder and Zhang \cite{balcan07} proved a
similar result for linear separators in higher dimensions,
and Castro and Nowak \cite{castro08} showed related improvements for
the space of boundary
fragment classes (under a somewhat stronger assumption than Tsybakov's).
However, these works left open the question of
whether such improvements could be achieved by an algorithm that does
not explicitly depend on the noise conditions (i.e., in the
\textit{agnostic} setting), and whether this type of improvement is achievable
for more general families of hypothesis classes, under the usual
complexity restrictions
(e.g., VC class, entropy conditions, etc.).
In a personal communication, John Langford and Rui Castro claimed $A^2$ achieves
these improvements for the special case of threshold classifiers
(a special case of this also appeared in \cite{balcan09}).
However, there remained an open question of whether such rate
improvements could be generalized to hold for arbitrary hypothesis
classes. In Section \ref{sec:rates}, we provide this
generalization. We analyze the rates achieved by $A^2$ under
Tsybakov's noise conditions \cite{mammen99,tsybakov04};
in particular, we find that these rates are strictly superior
to the known rates for passive learning, when the disagreement
coefficient is finite. We also study a novel modification of the algorithm
of Dasgupta, Hsu and Monteleoni \cite{dasgupta07}, proving that it
improves upon the rates of $A^2$
in its dependence on the disagreement coefficient.

Additionally, in Section \ref{sec:aggregation}, we address the second
limitation by proposing a general
model selection procedure for active learning with an arbitrary structure
of nested hypothesis classes. If the classes have restricted
expressiveness (e.g., VC classes),
the error rate for this algorithm converges to the best achievable error
by any classifier in the structure, at a rate that adapts to the noise
conditions
and complexity of the optimal classifier.
In general, if the structure is constructed to include arbitrarily
good approximations to any classifier, the error converges to the
Bayes error rate in the limit. In particular, if the Bayes optimal
classifier is in some class within the structure, the algorithm
performs nearly as well as running an agnostic active learning
algorithm on that single hypothesis class, thus preserving the convergence
rate improvements achievable for that class.

\section{Definitions and notation}
\label{sec:definitions}
In the active learning setting, there is an \textit{instance space}
$\mathcal{X}$,
a \textit{label space} $\Y= \{-1,+1\}$ and
some fixed distribution $\D_{XY}$ over $\mathcal{X} \times\Y$, with
marginal $\D_{X}$ over $\X$.
The restriction to binary classification ($\Y= \{-1,+1\}$) is intended
to simplify the discussion; however,
everything below generalizes quite naturally to multiclass
classification (where $\Y= \{1,2,\ldots,k\}$).

There are two sequences of random variables: $X_1, X_2, \ldots$ and
$Y_1, Y_2, \ldots,$ where each
$(X_i,Y_i)$ pair is independent of the others, and has joint
distribution $\D_{XY}$.
However, the learning algorithm is only permitted direct access to the
$X_i$ values (unlabeled data points),
and must request the $Y_i$ values one at a time, sequentially. That is,
the algorithm picks some index $i$ to observe the $Y_i$ value,
then after observing it, picks another index $i^\prime$ to observe the
$Y_{i^\prime}$ label value, etc.
We are interested in studying the rate of convergence of the error rate
of the classifier output by the learning
algorithm, in terms of the number of label requests it has made. %
%will also discuss convergence rates in terms of the number of
%unlabeled points the algorithm has observed.}}\ignore{
%we will specifically be interested in situations where that error rate
%converges to the best error rate achievable
%by any classifier, in the limit.}
To simplify the discussion, we will think of the data sequence as being
essentially inexhaustible,
and will study $(1-\delta)$-confidence bounds on the error rate of the
classifier produced by an algorithm permitted to
make at most $n$ label requests, for a fixed value $\delta\in
(0,1/2)$. The actual number of (unlabeled) data points the algorithm
uses will be made clear in the proofs (typically close to the number of
points needed by passive learning to achieve the
stated error guarantee).
%%% should these caveats, etc be moved elsewhere? is this too soon to
%get into minutia?

%Define $\Data_n = \{(X_1,Y_1),(X_2,Y_2),\ldots,(X_n,Y_n)\}$, a finite
%sequence consisting of the first $n$ examples.
%and define $\DataX_{n}^{\infty} = \{X_n,X_{n+1},\ldots\}$, a tail
%subsequence of the $X$ (unlabeled) samples.
A \textit{hypothesis class} $\C$ is any set of measurable
classifiers $h \dvtx\mathcal{X} \rightarrow\Y$. We will denote by
$d$ the VC dimension of $\C$ (see, e.g., \cite
{vapnik82,vapnik98,devroye96,blumer89,vapnik71}).
%, defined as the largest integer $m$ for which
%$\exists(x_1,x_2,\ldots,x_m) \in\X^m \text{ s.t. } |\{(y_1,y_2,
%= y_i\}| = 2^m$.
For any measurable $h \dvtx\X\rightarrow\Y$ and distribution $\D$ over
$\X\times\Y$,
define the \textit{error rate} of $h$ as $\er_{\D}(h) = \PP_{(X,Y)\sim
\D}\{h(X) \neq Y\}$; when $\D= \D_{XY}$,
we abbreviate this as $\er(h)$. This simply represents the risk under
the $0$--$1$ loss.
We also define the \textit{conditional error rate}, given a set $R
\subseteq\X$, as
$\er(h|R) = \PP\{h(X) \neq Y | X \in R\}$.
Let $\nu= \inf_{h \in\C} \er(h)$, called the \textit{noise rate} of
$\C$.
For any $x \in\mathcal{X}$, let
$\eta(x) = \PP\{Y=1 | X = x\}$, let $h^*(x) = 2 \mathbh{1}[\eta(x)
\geq1/2] - 1$ and let $\nu^* = \er(h^*)$.
We call $h^*$ the \textit{Bayes optimal classifier} and $\nu^*$ the
\textit{Bayes error rate}.
Additionally, define the \textit{diameter} of any set of classifiers $V$
as $\diam(V) = \sup_{h_1,h_2 \in V} \PP\{h_1(X) \neq h_2(X)\}$,
and for any $\varepsilon> 0$, define the diameter of the \textit
{$\varepsilon
$-minimal set} of $V$ as
$\diam(\varepsilon;V) = \diam(\{h \in V \dvtx \er(h) - \inf_{h^\prime\in
V}\er(h^\prime) \leq\varepsilon\})$.

For a classifier $h$, and a sequence $S = \{(x_1,y_1),(x_2,y_2),\ldots
,(x_m,y_m)\} \in(\X\times\Y)^m$, let
$\er_S(h) = \frac{1}{|S|}\sum_{(x,y) \in S} \mathbh{1}[h(x) \neq y]$
denote the \textit{empirical error rate} on $S$, [and define $\er_{\{\}
}(h) = 0$ by convention].
It will often be convenient to make use of sets of (index, label)
pairs, where the index is used to uniquely refer to an element of the
$\{X_i\}$ sequence
(while conveniently also keeping track of relative ordering information);
in such contexts, we will overload notation as follows.
For a classifier~$h$, and a finite set of (index, label) pairs $S = \{
(i_1, y_{1}),(i_2, y_{2}),\ldots,(i_m, y_{m})\} \subset\nats\times
\Y$,
let $\er_S(h) = \frac{1}{|S|} \sum_{(i,y) \in S} \mathbh
{1}[h(X_i) \neq y]$, (and $\er_{\{\}}(h) = 0$, as before).
Thus, $\er_S(h) = \er_{S^\prime}(h)$, where $S^\prime= \{(X_i,
y)\}_{(i,y) \in S}$.
For the indexed \textit{true} label sequence, $\Data^{(m)} = \{
(1,Y_1),(2,Y_2),\ldots,(m,Y_m)\}$,
we abbreviate this $\er_m(h) = \er_{\Data^{(m)}}(h)$, the
empirical error on the first $m$ data points.

In addition to the independent interest of understanding the rates
achievable here,
another primary interest in this setting is to quantify the achievable
\textit{improvements}, compared to \textit{passive learning}.
In this context, a passive learning algorithm can be formally defined
as a function
mapping the sequence $\{(X_1,Y_1), (X_2,Y_2),\ldots, (X_n,Y_n)\}$
to a classifier $\hat{h}_n$;
%in both applied and theoretical contexts;}% end ignore
for instance, perhaps the most widely studied family of passive
learning methods is that of \textit{empirical risk minimization} (e.g.,
\cite{vapnik82,vapnik98,koltchinskii06,massart06}),
which return a classifier $\hat{h}_n \in\arg\min_{h \in\C}
\er_n(h)$. For the purpose of this comparison, we
review known results on passive learning in several contexts below.

%s2.1 ###
\subsection{Tsybakov's noise conditions}

Here we describe a particular parametrization of noise distributions,
relative to a hypothesis class, often referred to as Tsybakov's noise
conditions \cite{mammen99,tsybakov04}, or margin conditions. These
noise conditions have recently received substantial attention in the
passive learning
literature, as they describe situations in which the asymptotic
minimax convergence rate of passive learning is faster than the worst
case $n^{-1/2}$
rate (e.g., \cite{mammen99,tsybakov04,koltchinskii06,massart06}).
\begin{condition}
\label{con:tsybakov}
There exist finite constants $\mu> 0$ and $\kappa\geq1$,
s.t.
$\forall\varepsilon>0$, $\diam(\varepsilon;\C) \leq\mu\varepsilon^{
{1/\kappa}}$.
\end{condition}

This condition is satisfied when, for example,
\[
\exists\mu^{\prime} > 0, \kappa\geq1 \mbox{ s.t. } \exists h \in
\C\dvtx\forall h^\prime\in\C\qquad \er(h^\prime) - \nu\geq\mu^{\prime}
\PP\{h(X) \neq h^\prime(X)\}^{\kappa},
\]
\cite{koltchinskii06}. It is also satisfied when the Bayes optimal
classifier is in $\C$ and
\[
\exists\mu^{\prime\prime} > 0, \alpha\in(0,\infty) \mbox{ s.t. }
\forall\varepsilon> 0\qquad \PP\{ |\eta(X) - 1/2| \leq\varepsilon\} \leq\mu
^{\prime\prime} \varepsilon^{\alpha},
\]
where $\kappa$ and $\mu$ are functions of $\alpha$ and $\mu^{\prime
\prime}$ \cite{mammen99,tsybakov04};
in particular, $\kappa= (1+\alpha)/\alpha$.
As we will see, the case where $\kappa= 1$ is particularly
interesting; for instance, this is the case when $h^* \in\C$ and
$\PP\{|\eta(X)-1/2| > c\}=1$ for some constant $c\in(0,1/2)$.
Informally, in many cases Condition \ref{con:tsybakov} can be realized in
terms of the relation between magnitude of noise and distance to the
optimal decision boundary; that is, since in practice the amount of
noise in
a data point's label is often inversely related to the distance from the
decision boundary, a small $\kappa$ value may often result from
having low density near the decision boundary (i.e., large margin);
when this is not the case, the value of $\kappa$ is often
determined by how quickly $\eta(x)$ changes as $x$ approaches the
decision boundary. See
\cite{mammen99,tsybakov04,koltchinskii06,castro08,massart06,balcan07}
for further interpretations of this condition.

It is known that when this condition is satisfied for some $\kappa
\geq1$ and $\mu> 0$, the passive learning method of empirical risk
minimization achieves a convergence rate guarantee, holding with
probability $\geq1-\delta$, of
\[
\er\Bigl(\mathop{\arg\min}_{h \in\C} \er_n(h)\Bigr) - \nu\leq c \biggl(\frac{d
\log n +
\log(1/\delta)}{n}\biggr)^{{\kappa}/({2\kappa- 1})},
% should I include dependence on $\mu$?
\]
where $c$ is a ($\kappa$ and $\mu$-dependent)
constant (this follows from \cite{koltchinskii06,massart06};
see Appendix B of the supplementary material \cite{hanneke10bsup},
especially (17) and Lemma 5, for the details). Furthermore,
for some hypothesis classes, this is known to be a tight bound (up to
the log factor) on the minimax convergence rate, so that there is
\textit{no} passive learning algorithm for these classes for which we
can guarantee a faster convergence rate, given that the guarantee
depends on $\D_{XY}$ only through $\mu$ and
$\kappa$ \cite{castro08,tsybakov04} (see also Appendix D of
the supplementary material \cite{hanneke10bsup}).

%s2.2 ###
\subsection{Disagreement coefficient}
\label{subsec:disagreement-coefficient}

The disagreement coefficient, introduced in \cite{hanneke07b}, is a
measure of
the complexity of an active learning problem, which has proven quite
useful for
analyzing the convergence rates of certain types of active learning algorithms:
for example, the algorithms of \cite
{cohn94,balcan06,dasgupta07,beygelzimer09}.
Informally, it quantifies how much disagreement there is among a set of
classifiers
relative to how close to some $h$ they are.
The following is a version of its definition, which we will use
extensively below.
For any hypothesis class $\C$ and $V \subseteq\C$, let
\[
\operatorname{DIS}(V) = \{x \in\X\dvtx\exists h_1,h_2 \in V \mbox{ s.t. }
h_1(x) \neq h_2(x)\}.
\]
For $r \in[0,1]$ and measurable $h \dvtx\X\rightarrow\Y$, let
\[
B(h,r) = \bigl\{h^\prime\in\C\dvtx\PP\{h(X) \neq h^\prime(X)\} \leq r\bigr\}.
\]
%
%when $\C$ or $\D$ is clear from the context, we may simply write $B_{
%Let $\operatorname{DIS}(B(h,r)) = \{x \in\X: \exists h^\prime\in
%
\begin{definition}
\label{def:disagreement-coefficient}
The \textit{disagreement coefficient} of $h$ with respect to $\C$ under
$\D_{X}$
is defined as
\[
\theta_{h} = \sup_{r > r_0} \frac{\PP(\operatorname{DIS}(B(h,r)))}{r},
\]
where $r_0 = 0$ (though see Appendix \ref{subsec:r0} for alternative
possibilities for $r_0$).
\end{definition}
\begin{definition}
\label{def:global-disagreement-coefficient}
We further define the disagreement coefficient for the hypothesis
class $\C$ with respect to the target distribution $\D_{XY}$ as
$\theta= \liminf_{k \rightarrow\infty} \theta_{h^{[k]}}$, where
$\{h^{[k]}\}$ is any sequence in $\C$ with $\er(h^{[k]})$ monotonically
decreasing to $\nu$; [by convention, take every $h^{[k]}\in\arg\min
_{h \in\C} \er(h)$ if the minimum is achieved].
\end{definition}

In Definition \ref{def:disagreement-coefficient}, it is conceivable
that $\operatorname{DIS}(B(h,r))$ may sometimes not be measurable.
In such cases, we can define $\PP(\operatorname{DIS}(B(h,r)))$ as the \textit{outer}
measure \cite{vanderVaart96}, so that
it remains well defined. We continue this practice below, letting $\PP
$ and $\E$ (and indeed any
reference to ``probability'') refer to the outer expectation and
measure in any context for which this is necessary.\looseness=1

Because of its simple intuitive interpretation, measuring the amount
of disagreement in a local neighborhood of some classifier $h$, the
disagreement coefficient has the wonderful property of being
relatively simple to calculate for a wide range of learning problems,
especially when those problems have a natural geometric
representation. To illustrate this, we will go through a few simple
examples from \cite{hanneke07b}.

Consider the hypothesis class of thresholds $h_z$ on the interval $[0,1]$
[for $z \in(0,1)$], where $h_z(x) = +1$ iff $x \geq z$.
Furthermore, suppose $\D_{X}$ is uniform on $[0,1]$.
In this case, it is clear that the disagreement coefficient is $2$,
since for sufficiently small $r$, the region of disagreement of
$B(h_z,r)$ is
$[z-r,z+r)$, which has probability mass $2r$.
In other words, since the disagreement region grows with $r$
in two disjoint directions, each at rate $1$, we have $\theta_{h_z} = 2$.

As a second example, consider the disagreement coefficient for \textit
{intervals}
on $[0,1]$. As before, let $\X= [0,1]$ and $\D_{X}$ be uniform, but
this time $\C$ is the set of intervals $h_{[a,b]}$ such that
for $x \in[0,1]$, $h_{[a,b]}(x) = +1$ iff $x \in[a,b]$ (for $0 < a <
b < 1$).
In contrast to thresholds, the disagreement coefficients $\theta
_{h_{[a,b]}}$ for the space of intervals vary widely depending on the
particular $h_{[a,b]}$.
Specifically, we have
$\theta_{h_{[a,b]}} = \max\{\frac{1}{b-a}, 4\}$.
To see this, note that when $0 < r < b-a$, every interval in
$B(h_{[a,b]},r)$ has its lower and upper boundaries within $r$ of $a$
and $b$, respectively;
thus, $\PP(\operatorname{DIS}(B(h_{[a,b]},r))) \leq4 r$, with equality for
sufficiently small $r$. However, when $r > b-a$, \textit{every} interval
of $\mathrm{width} \leq r-(b-a)$ is in $B(h_{[a,b]},r)$,
so that $\PP(\operatorname{DIS}(B(h_{[a,b]},r))) = 1$.

As a slightly more involved example, \cite{hanneke07b} studies the
scenario where $\X$ is the surface of the origin-centered unit sphere
in $\mathbb{R}^d$ for $d > 2$, $\C$ is the space of all linear separators
whose decision surface passes through the origin, and $\D_X$ is the uniform
distribution on $\X$; in this case, it turns out $\forall h \in\C$ the
disagreement coefficient $\theta_h$ satisfies
\[
\frac{\pi}{4} \sqrt{d} \leq\theta_h \leq\pi\sqrt{d}.
\]
%
%Wortman \cite{hanneke08a} study the
%scenario where $\X$ is the surface of the origin-centered unit sphere
%in $\mathbb{R}^d$ for $d \geq2$, $\C$ is the space of all linear
%separators and $\D_X$ is the uniform distribution on $\X$; in this
%case,
%it turns out $\forall h \in\C$ the disagreement coefficient $
%%c_1 \max\{\frac{1}{p_h}, \sqrt{d}\} \leq\theta_h \leq c_2
%where $p_h = \min\{\PP(h(X) = +1),\PP(h(X) = -1)\}$ is the minority
%class probability, %and $c_1,c_2 > 0$ are universal constants.
%and $c > 0$ is a universal constant.}% end ignore

The disagreement coefficient has many interesting properties that can
help to bound its value for a given hypothesis class and distribution.
We list a few elementary properties below. Their proofs, which are
quite short and follow directly from the definition, are left as easy
exercises.
\begin{lemma}[(Close marginals \cite{hanneke07b})]
\label{lem:close}
 Suppose $\exists\lambda\in
(0,1]$ s.t. for any measurable set $A \subseteq\X$, $\lambda
\PP_{\D_{X}}(A) \leq\PP_{\D_X^\prime}(A) \leq
\frac{1}{\lambda}\PP_{\D_X}(A)$. Let $h \dvtx\X\rightarrow\Y$ be
a measurable classifier, and suppose $\theta_h$ and $\theta_h^\prime$
are the disagreement coefficients for $h$ with respect to $\C$ under
$\D_X$ and $\D_X^\prime$, respectively. Then
%$\lambda\PP_{\D_X}(\operatorname{DIS}(\dense{B}_{\D_X}(h,\lambda r))) \leq\PP_{\D_X^
%and thus
%
\[
\lambda^2 \theta_h \leq\theta_h^\prime\leq\frac{1}{\lambda
^2}\theta_h.
\]
\end{lemma}
%
%For any $h^\prime\in\C$,
%Thus, $B_{\D_X}(h,\lambda r) \subseteq B_{\D_X^\prime}(h,r) \subseteq
%B_{\D_X}(h,r/\lambda)$, which implies
%The result on the disagreement coefficients then follows from the
%definition.
%journal version)
%
\begin{lemma}[(Finite mixtures)]\label{lem:mixtures}
 Suppose $\exists\alpha\in[0,1]$ s.t. for any
measurable set $A \subseteq\X$,
$\PP_{\D_X}(A) = \alpha\PP_{\D_1}(A) + (1-\alpha)\PP_{\D_2}(A)$.
For a measurable $h \dvtx\X\rightarrow\Y$, let
$\theta_h^{(1)}$ be the disagreement coefficient with respect to $\C$
under $\D_1$,
$\theta_h^{(2)}$ be the disagreement coefficient with respect to $\C$
under $\D_2$,
and $\theta_h$ be the disagreement coefficient with respect to $\C$
under $\D_{X}$. Then
\[
\theta_h \leq\theta_h^{(1)} + \theta_h^{(2)}.
\]
%
%%% TODO: can I do a lower bound \min\{\theta_h^{(1)},\theta_h^{(2)}\}
%???
\end{lemma}
%
%%For any $r > 0$, \\
%$\PP_{\D}\{\operatorname{DIS}(\dense{B}_{\D}(h,r))\}$
%& = \alpha\mathbb{P}_{\D_1}\{\operatorname{DIS}(\dense{B}_{\D}(h,r))\} + (1- \alpha)
%& \leq\alpha\mathbb{P}_{\D_1}\{\operatorname{DIS}(\dense{B}_{\D_1}(h,r/\alpha))\} +
%(1- \alpha) \mathbb{P}_{\D_2}\{\operatorname{DIS}(\dense{B}_{\D_2}(h,r/(1-\alpha)))\}
%& \leq\theta_h^{(1)}r + \theta_h^{(2)}r. %\qedhere
%%The result then follows from the definition of disagreement
%coefficient.
%journal version)
%
\begin{lemma}[(Finite unions)]\label{lem:unions}
Suppose $h \in\C_1 \cap\C_2$ is a classifier
s.t. the disagreement coefficient with respect to $\C_1$ under $\D_X$ is
$\theta_h^{(1)}$ and with respect to $\C_2$ under $\D_X$ is
$\theta_h^{(2)}$. Then if $\theta_h$ is the disagreement coefficient
with respect to $\C= \C_1 \cup\C_2$ under $\D_X$, we have that
\[
\max\bigl\{\theta_h^{(1)},\theta_h^{(2)}\bigr\} \leq\theta_h
\leq\theta_h^{(1)} + \theta_h^{(2)}.
\]
In fact, even if $h \notin\C_1 \cap\C_2$, we still have $\theta_h
\leq\theta_h^{(1)} + \theta_h^{(2)} + 2$.
\end{lemma}
%
%The lower bound is clear from the definition of disagreement
%coefficient.
%For the upper bound, we must have that $\operatorname{DIS}(\dense{B}_{\C_1 \cup
%and thus
%journal version)
\eject

See \cite
{hanneke07b,dasgupta07,hanneke08a,beygelzimer09,friedman09,wang09}
for further discussions of various uses of the disagreement
coefficient and related notions and extensions in active learning. In
particular, Friedman \cite{friedman09} proves that any hypothesis class
and distribution satisfying certain general regularity conditions
will admit finite constant bounds on $\theta$. Also, Wang \cite{wang09}
bounds the disagreement coefficient for certain nonparametric
hypothesis classes, characterized by smoothness of their decision surfaces.
Additionally, Beygelzimer, Dasgupta and Langford \cite{beygelzimer09}
present an interesting analysis
using a natural extension of the disagreement coefficient to study
active learning with a larger family of loss functions beyond $0$--$1$
loss. %\ignore{As a related aside, although the focus of this paper is
%active
%learning, interestingly the disagreement coefficient also has
%applications in the analysis of \textit{passive} learning; see
%(I'm not going to publish that one before submitting this)

%%% TODO: brief discussion of suboptimality of disagreement coefficient
%in terms of constants ?

The disagreement coefficient has deep connections to several other
quantities, such as doubling dimension \cite{long07} and VC
dimension \cite{vapnik82}.
Additionally, a related quantity, referred to as the ``capacity function,''
was studied in the 1980s by Alexander in the passive learning literature,
in the context of ratio-type empirical processes \cite
{alexander85,alexander86,alexander87}
and recently was further developed
by Gin\'{e} and Koltchinskii \cite{gine06}; interestingly, in this
latter work, Gin\'{e} and Koltchinskii study
a localized version of the capacity function,
which in our present context can essentially be viewed as the function
$\tau(r) = \PP(\operatorname{DIS}(B(h,r)))/r$, so that $\theta_h = \sup_{r > r_0}
\tau(r)$.

%s3 ###
\section{General algorithms}
\label{sec:algorithms}

We begin the discussion of the algorithms we will analyze by noting
the underlying inspiration that unifies them. Specifically, at this
writing, all of the published general-purpose agnostic active learning
algorithms achieving nontrivial improvements are derivatives of a
basic technique proposed by Cohn, Atlas and Ladner \cite{cohn94} for
the realizable active
learning problem. Under the assumption that there exists a perfect
classifier in $\C$, they proposed an algorithm which processes
unlabeled data points in sequence, and for each one it determines whether
there is a classifier in $\C$ consistent with all previously
observed labels that predicts $+$1 for this new point
\textit{and} one that predicts $-$1 for this new point;
if so, the algorithm requests the label, and
otherwise it does not request the label; after $n$ label requests, the
algorithm returns any classifier consistent with all observed labels.
In some sense, this algorithm corresponds to the very least we could
expect of an active learning algorithm, as it never requests the label
of a point it can derive from known information, but otherwise
makes no effort to search for informative data points. The idea is appealing,
not only for its simplicity, but also for its extremely
efficient use of unlabeled data; in fact, under the stated assumption,
the algorithm produces a classifier consistent with the labels of
\textit{all}
of the unlabeled data it processes, including those it does \textit{not}
request the labels of.

We can equivalently think of this algorithm as maintaining two sets: $V
\subseteq\C$ is the set of candidate hypotheses still under
consideration, and $R = \operatorname{DIS}(V)$ is their region of disagreement. We
can then think of the algorithm as requesting a random labeled point
from the conditional distribution of $\D_{XY}$ given that $X \in R$,
and subsequently removing from $V$ any classifier inconsistent with
the observed label. A formal definition of the algorithm
is given as follows.

\vspace*{24pt}

\begin{bigboxit}
\textbf{Algorithm 0} \\
Input: hypothesis class $\C$, label budget $n$\\
Output: classifier $\hat{h}_n \in \C$\\
%{\vskip -2mm}\line(1,0){377}\\
{\vskip -2mm}\line(1,0){346}\\
0. $V_0 \leftarrow \C$, $t \leftarrow 0$\\
1. For $m = 1,2,\ldots$\\
%2. \quad If $\exists h_1,h_2 \in V_t$ s.t. $h_1(X_m) \neq h_2(X_m)$,\\
2. \quad If $X_m \in \operatorname{DIS}(V_t)$,\\
3. \quad\quad Request $Y_m$\\
4. \quad\quad $t \leftarrow t+1$\\
5. \quad\quad $V_t \leftarrow \{h \in V_{t-1} : h(X_m) = Y_m\}$\\
6. \quad If $t = n$ or $\{m^\prime > m : X_{m^\prime} \in \operatorname{DIS}(V_t)\}=\varnothing$, Return any $\hat{h}_n \in V_t$
%
%0. $V_0 \leftarrow \C$, $t \leftarrow 0$, $m \leftarrow 0$ \\
%1. For $t = 1,2,\ldots,n$\\
%2. \quad If $t=n+1$ or $\{m^\prime > m  : X_{m^\prime} \in DIS(V_{t-1})\} = \varnothing$, \\
%3. \quad\quad Return an arbitrary $\hat{h}_n \in V_{t-1}$\\
%4. \quad $m \leftarrow \min\{m^\prime > m : X_{m^\prime} \in DIS(V_{t-1})\}$\\
%5. \quad Request $Y_m$\\
%6. \quad $V_t \leftarrow \{h \in V_{t-1} : h(X_m) = Y_m\}$\\
\end{bigboxit}

\vspace*{24pt}

The algorithms described below for the problem of active learning with
label noise each represent noise-robust variants of this basic idea.
They work to reduce the set of candidate hypotheses, while only
requesting the labels of points in the region of disagreement of
these candidates. The trick is to only remove a classifier from the
candidate set once we have high statistical confidence that it is
worse than some other candidate classifier so that we never remove the
best classifier. However, the two algorithms differ somewhat in the
details of how that confidence is calculated.

%during each of which the sampling distribution is fixed to the
%conditional distribution given the region of disagreement at the
%start of the phase; it samples labeled points according to
%this conditional distribution, and uses them to calculate upper
%and lower bounds on the error rates of the candidate classifiers
%under that distribution; the algorithm removes a classifier from
%consideration if its lower bound exceeds another classifier's upper
%bound, indicating the former is inferior to the latter; it transitions
%between phases when the region of disagreement of the candidate set
%has reduced significantly in size, enabling more focused sampling in
%the next phase.
%
%%%% TODO: this really isn't a very good description of the algorithm
%%%% but I'll leave it like this for now, since I'll probably have to
%cut it
%%%% for space anyway.
%%%% If there is space, I can probably lift a good description from my
%%%% thesis proposal
%The second algorithm maintains the candidate set of classifiers
%algorithm). It processes each data point in sequence, and using a
%confidence bound designed specifically for this situation, determines
%whether the point is in the region of disagreement of the remaining
%set of candidate hypotheses; however, in contrast to the first
%algorithm, the confidence bound is calculated in terms of the total
%number of data points processed so far, including those it did not
%request the labels of.}% end ignore

%can potentially update the sampling distribution after every label
%request, and %calculates bounds in terms of the total number of

%s3.1 ###
\subsection{Algorithm 1}

The first noise-robust algorithm we study, originally proposed by
Balcan, Beygelzimer and Langford \cite{balcan06}, is
typically referred to as $A^2$ for \textit{Agnostic Active}. This was
historically the first general-purpose agnostic active learning
algorithm shown to achieve improved error guarantees for certain
learning problems in certain ranges of $n$ and $\nu$.
Below is a variant of this algorithm.
It is defined in terms of two functions: $\UB $ and $\LB $.
These represent upper and lower confidence bounds on the error rate of
a classifier from $\C$ with respect to an arbitrary sampling
distribution, as a function of a labeled sequence sampled according to
that distribution.
Some steps in the algorithm require calculating certain probabilities,
such as $\PP(\operatorname{DIS}(V))$ or $\PP(R)$; later, we discuss replacing these with
appropriate estimators.\vspace*{4pt}
\vfill

\begin{bigboxit}
\textbf{Algorithm 1}\\
Input: hypothesis class $\C$, label budget $n$, confidence $\delta$, functions $\UB$ and $\LB$\\
Output: classifier $\hat{h}_n$\\
{\vskip -2mm}\line(1,0){346}\\
%{\vskip -2mm}\line(1,0){377}\\
%{\vskip -1mm}
0. $V \leftarrow \C$, $R \leftarrow \operatorname{DIS}(\C)$, $Q \leftarrow \varnothing$, $m \leftarrow 0$\\
1. For $t = 1,2,\ldots,n$\\
2. \quad If $\PP(\operatorname{DIS}(V)) \leq \frac{1}{2}\PP(R)$\\
3. \quad\quad $R \leftarrow \operatorname{DIS}(V)$; $Q \leftarrow \varnothing$\\
4. \quad\quad If $\PP(R) \leq 2^{-n}$, Return any $\hat{h}_n \in V$\\
5. \quad $m \leftarrow \min\{m^\prime > m \dvtx X_{m^\prime} \in R\}$\\
6. \quad Request $Y_m$ and let $Q \leftarrow Q \cup \{(m,Y_m)\}$\\
7. \quad $V \leftarrow \{h \in V \dvtx \LB(h,Q,\delta/n) \leq \min_{h^\prime \in V}\UB(h^\prime,Q,\delta/n)\}$\\
%8. \quad If $\left(\min_{h \in V} \UB(h,Q,\delta/n) - \min_{h \in V}LB(h,Q,\delta/n)\right)\mathbb{P}(R) < \beta$\\
8. \quad $h_t \leftarrow \arg\min_{h \in V} \UB(h,Q,\delta/n)$\\
9. \quad $\beta_t \leftarrow (\UB(h_t,Q,\delta/n) - \min_{h \in V}\LB(h,Q,\delta/n))\mathbb{P}(R)$\\
10. Return $\hat{h}_n = h_{\hat{t}}$, where $\hat{t} = \arg\min_{t \in \{1,2,\ldots,n\}} \beta_t$
\end{bigboxit}

\vspace*{6pt}

The intuitive motivation behind the algorithm is the following.
It focuses on reducing the set of candidate hypotheses $V$, while being
careful not to throw away the best classifier $h_{\C}^{*} = \arg\min
_{h \in\C} \er(h)$
(supposing, for this informal explanation, that $h_{\C}^{*}$ exists).
Given that this is satisfied at any given time in the algorithm, it
makes sense to focus our samples to the region $\operatorname{DIS}(V)$,
since a classifier $h_1 \in V$ has smaller error rate than another
classifier $h_2 \in V$ if and only if it has smaller conditional error
rate given $\operatorname{DIS}(V)$.
For this reason, on each round, we seek to remove from $V$ any $h$ for
which our confidence bounds indicate that $\er(h | \operatorname{DIS}(V)) > \er(h_{\C
}^{*} | \operatorname{DIS}(V))$.
However, so that we can make use of known results for i.i.d. samples,
we freeze the sampling region $R \supseteq \operatorname{DIS}(V)$ and collect an i.i.d.
sample from the conditional given this region, updating the region only
when doing so allows us to further significantly focus the samples;
for this same reason, we also reset the collection of samples $Q$ every
time we update the region $R$, so that it represents samples from the
conditional given $R$. Finally, we maintain the values $\beta_t$,
which represent confidence upper bounds on $\er(h_t) - \nu= (\er(h_t |
R) - \er(h_{\C}^{*} | R)) \PP(R)$,
and we return the $h_t$ minimizing this confidence bound; note that it
does not suffice to return $h_n$, since the final $Q$ set might be small.

As long as the confidence bounds $\UB $ and $\LB $ satisfy (overloading notation
in the natural way)
%in the natural way, so that $\UB (h,\{(i_j,Y_{i_j})\}_j,\delta^\prime) =
%footnote so I can just put it in parens)
%
\[
\PP_{Z\sim\D^m}\{\forall h \in\C, \LB (h,Z,\delta^\prime) \leq
\er_{\D}(h) \leq \UB (h,Z,\delta^\prime)\} \geq1-\delta^\prime
\]
for any distribution $\D$ over $\X\times\Y$ and any $\delta^\prime
\in
(0,1)$, and $\UB $ and $\LB $ converge to each other as $m$ grows, it is
known that
a $1-\delta$ confidence bound on $\er(\hat{h}_n) - \nu$
converges to $0$ \cite{balcan06}. For instance,
Balcan, Beygelzimer and Langford \cite{balcan06} suggest defining
these functions based on classic
results on uniform convergence rates in passive
learning \cite{vapnik82}, such as
%
%e1 ###
%
\begin{eqnarray}
\label{eqn:bbl-bounds}
\UB (h,Q,\delta^\prime) &=& \min\{\er_Q(h) + G(|Q|,\delta^\prime),1\}
,\nonumber\\[-8pt]\\[-8pt]
\LB (h,Q,\delta^\prime) &=& \max\{\er_Q(h) - G(|Q|,\delta^\prime),0\}
,\nonumber
\end{eqnarray}
where $G(m,\delta^\prime) = \frac{1}{m} + \sqrt{\frac{\ln
({4}/{\delta^\prime}) + d \ln({2 e m}/{d})}{m}}$ for $m \geq d$,
and by convention $G(m,\delta^\prime) = \infty$ for $m < d$.
This choice of $\UB $ and $\LB $ is motivated by the following lemma, due
to Vapnik \cite{vapnik98}.
\begin{lemma}
\label{lem:uniform}
For any distribution $\D$ over $\X\times\Y$, and any $\delta
^\prime\in(0,1)$ and $m \in\nats$,
with probability $\geq1-\delta^\prime$ over the draw of $Z \sim\D^m$,
every $h \in\C$ satisfies
%
%e2 ###
%
\begin{equation}\label{eqn:uniform}
|\er_{Z}(h) - \er_{\D}(h)| \leq G(m,\delta^\prime).
%m}{d}}{m}}.
%%% TODO: in the journal version, why not use this tighter bound? I
%need to make both the defn of \UB /\LB  and this lemma match, without
%making it all look weird.
\end{equation}
\end{lemma}

To avoid computational issues, instead of explicitly representing the
sets $V$ and~$R$, we may implicitly represent them as a set of
constraints imposed by the condition in step 7 of previous iterations.
We may also replace $\PP(\operatorname{DIS}(V))$ and $\PP(R)$ by estimates, since these
quantities can be estimated to arbitrary precision with arbitrarily
high confidence using only \textit{unlabeled} data. Specifically,
the convergence rates proven below can be preserved up to constant factors
by replacing these quantities with confidence bounds based on a finite number
of unlabeled data points%\ignore{ polynomial in the reciprocals of the
%stated excess
%error guarantee and $\delta$}
; the details of this are included in
Appendix C of the supplementary material~\cite{hanneke10bsup}.
As for the number of unlabeled data points required by the above
algorithm itself,
note that if $\PP(\operatorname{DIS}(V))$ becomes small, it will use a large number
of unlabeled data points;
however, $\PP(\operatorname{DIS}(V))$ being small also indicates $\er(\hat{h}_n) -
\nu$ is small (and indeed $\beta_t$).
In particular, to get an excess error rate of $\varepsilon$, the
algorithm will generally require a number of unlabeled data points only
polynomial in $1/\varepsilon$;
also, the condition in step 4 guarantees the total number of unlabeled
data points used by the algorithm is bounded with high probability.
For comparison, recall that passive learning typically requires a
number of \textit{labeled} data points polynomial in $1/\varepsilon$.

%s3.2 ###
\subsection{Algorithm 2}

The second noise-robust algorithm we study was originally proposed by
Dasgupta, Hsu and Monteleoni \cite{dasgupta07}. It uses a type of
constrained passive learning subroutine, \textsc{Learn},
defined as follows for two sets of labeled data points, $\LL$ and $Q$.
\[
\mbox{\textsc{Learn}}_{\C}(\LL,Q) = \mathop{\arg\min}_{h \in\C\dvtx
\er_{\LL}(h)=0}\er_Q(h).
\]
By convention, if no $h \in\C$ has $\er_{\LL}(h)=0$, $\mbox{\textsc
{Learn}}_{\C}(\LL,Q)=\varnothing$.
The algorithm is formally defined below, in terms of a sequence of
estimators $\Delta_{m}$, defined later.%\vspace*{4pt}

\begin{bigboxit}
\textbf{Algorithm 2}\\
Input: hypothesis class $\C$, label budget $n$, confidence $\delta$,
functions $\Delta_{m}$\\
Output: classifier $\hat{h}_n$, sets of (index, label) pairs $\LL$
and $Q$\\
{\vskip-2mm}\line(1,0){346}\\
%{\vskip-2mm}\line(1,0){377}\\
%{\vskip-1mm}
0. $\LL\leftarrow\varnothing$, $Q \leftarrow\varnothing$\\
1. For $m = 1,2,\ldots$\\
2. \quad If $|Q|=n$ or $m > 2^n$, Return $\hat{h}_n = \mbox{\textsc
{Learn}}_{\C}(\LL,Q)$ along with $\LL$ and $Q$\\
3. \quad For each $y \in\{-1,+1\}$, let $h^{(y)} = \mbox{\textsc{Learn}}_{\C
}(\LL\cup\{(m,y)\},Q)$\\
4. \quad If some $y$ has $h^{(-y)} = \varnothing$ or \\
{\hskip2.5cm}$\er_{\LL\cup Q}\bigl(h^{(-y)}\bigr)-\er_{\LL\cup Q}\bigl(h^{(y)}\bigr) >
\Delta_{m-1}\bigl(\LL,Q,h^{(y)},h^{(-y)},\delta\bigr)$\\
5.\qquad  Then $\LL\leftarrow\LL\cup\{(m,y)\}$\\
6. \quad Else Request the label $Y_m$ and let $Q \leftarrow Q \cup\{
(m,Y_m)\}$
\end{bigboxit}

\vspace*{6pt}

The algorithm maintains two sets of labeled data points: $\LL$ and
$Q$. The set $Q$ represents
points of which we have requested the labels. The set $\LL$ represents
the remaining points, and
the labels of points in $\LL$ are \textit{inferred}. Specifically,
suppose (inductively) that at some time $m$
we have that every $(i,y) \in\LL$ has $h_{\C}^{*} (X_i) = y$, where
$h_{\C}^{*} = \arg\min_{h \in\C} \er(h)$
(supposing the min is achieved, for this informal motivation). At any
point, we can be
fairly confident that $h_{\C}^{*}$ will have relatively small
empirical error rate. Thus, if all of the classifiers
$h$ with $\er_{\LL}(h) = 0$ and $h(X_{m}) = -y$ have relatively large
empirical error rates compared to some
$h$ with $\er_{\LL}(h) = 0$ and $h(X_{m}) = y$, we can confidently
infer that $h_{\C}^{*}(X_{m}) = y$.
Note that this is not the \textit{true} label $Y_{m}$, but a sort of
``denoised'' version of it.
Once we infer this label, since we are already confident that this is
the $h_{\C}^{*}$
label, and $h_{\C}^{*}$ is the classifier we wish to compete with, we
simply add this label as a \textit{constraint}:
that is, we require every classifier under consideration in the future
to have $h(X_{m}) = h_{\C}^{*}(X_{m})$.
This is how elements of $\LL$ are added. On the other hand, if we
cannot confidently infer $h_{\C}^{*}(X_{m})$, because
some classifiers labeling $X_m$ as $-h_{\C}^{*}(X_{m})$ also have relatively small
empirical error rates, then we simply request the label $Y_{m}$
and add it to the set $Q$. Note that in order to make this comparison,
we needed to be able to calculate the differences
of empirical error rates; however, as long as we only consider the set
of classifiers $h$ that \textit{agree} on the labels in $\LL$,
we will have $\er_{\LL\cup Q}(h_1) - \er_{\LL\cup Q}(h_2) = \er_{m}(h_1)
- \er_{m}(h_2)$, for any two such classifiers $h_1$ and~$h_2$,
where $m = |\LL\cup Q|$.

The key to the above argument is carefully choosing a threshold for how
large the difference in empirical error rates needs to be before
we can confidently infer the label.
For this purpose, Algorithm 2 is defined in terms of a function,
$\Delta_m(\LL,Q,h^{(y)},h^{(-y)},\delta)$,
representing a threshold for a type of hypothesis test. This threshold
must be set
carefully, since the sequence of labeled data points corresponding to
$\LL\cup Q$ is not actually an i.i.d. sample from $\D_{XY}$.
Dasgupta, Hsu and Monteleoni \cite{dasgupta07} suggest defining this
function as
%
%e3 ###
%
\begin{equation}
\label{eqn:dhm-bound}\quad
\Delta_m\bigl(\LL,Q,h^{(y)},h^{(-y)},\delta\bigr) = \beta_m^2 + \beta_m
\bigl(\sqrt{\er_{\LL\cup Q}\bigl(h^{(y)}\bigr)}+\sqrt{\er_{\LL\cup
Q}\bigl(h^{(-y)}\bigr)}\bigr),
\end{equation}
where $\beta_m = \sqrt{\frac{4\ln(8m(m +
1)\mathcal{S}(\C,2m)^2/\delta)}{m}}$ and $\mathcal{S}(\C,2m)$ is the
shatter coefficient (e.g., \cite{vapnik98,devroye96}); this
suggestion is based on a confidence bound
they derive, and they prove the correctness of the algorithm with this
definition, meaning that the $1-\delta$ confidence bound on its error
rate converges to $\nu$ as $n \to\infty$.
For now we will focus on the first return value (the classifier),
leaving the others for Section \ref{sec:aggregation}, where they will
be useful for chaining multiple executions together.

%s4 ###
\section{Convergence rates}
\label{sec:rates}

In both of the above cases, one can prove guarantees stating
that neither algorithm's convergence rates are ever significantly worse
than passive
learning by empirical risk minimization \cite{balcan06,dasgupta07}.
However, it is even more
interesting to discuss situations in which one can prove error rate
guarantees for these algorithms significantly \textit{better} than those
achievable by passive learning. In this section, we begin by
reviewing known results on these potential improvements, stated in
terms of the disagreement coefficient; we then proceed to discuss new
results for Algorithm 1 and a novel variant of Algorithm 2, and
describe the convergence
rates achieved by these methods in terms of the disagreement coefficient
and Tsybakov's noise conditions.

To simplify the presentation,
for the remainder of this paper we will restrict the discussion to situations
with $\theta> 0$ (and therefore $\C$ with $d > 0$ too). Handling the
extra case of $\theta=0$
is a trivial matter, since $\theta= 0$ would imply that any proper learning
algorithm achieves excess error $0$ for all values of $n$.

%s4.1 ###
\subsection{The disagreement coefficient and active learning: Basic
results} %%% TODO: find a better section title

Before going into the results for general distributions $\D_{XY}$ on
$\X\times\Y$, it will be instructive to first look at the
special case when the noise rate is zero. Understanding how the
disagreement coefficient enters into the analysis of this simpler case
may aid in digestion of the theorems and proofs for the general case
presented later, where it plays an essentially analogous role. Most
of the major ingredients of the proofs for the general case can be
found in this special case, albeit in a much simpler form. Although
this result has not previously been published, the proof is
essentially analogous to (one case of) the analysis of
Algorithm 1 in \cite{hanneke07b}.
\begin{theorem}
\label{thm:cal-upper}
Let $f \in\C$ be such that $\er(f) = 0$ and $\theta_f < \infty$.
$\forall n \in\nats$ and $\delta\in(0,1)$, with probability $\geq
1-\delta$ over the draw of the unlabeled data,
the classifier $\hat{h}_n$ returned by Algorithm 0 after $n$ label
requests satisfies
\[
\er(\hat{h}_n) \leq2\cdot \exp\biggl\{-\frac{n}{12 \theta_f (d
\ln(22 \theta_f) + \ln(3n/\delta)
)}\biggr\}.
%%% TODO: I should be able to optimize constants here by considering
%different interval lengths; I should get the general formula for the
%bound and solve for the optimal value (not just 2)
\]
\end{theorem}
\begin{pf}%[Proof of Theorem \ref{thm:cal-upper}]
%The case $diam(\C) = 0$ is trivial, so assume $diam(\C) > 0$ (and thus
%$d \geq1$ and $\theta_f \geq1$).
%We will show that if $V_t$ is the set of classifiers in $\C$
%consistent with the first $t$ label requests,
%then $diam(V_n)$ is at most the above bound value, which implies the
%result.
As in the algorithm, let $V_t$ denote the set of classifiers in $\C$
consistent with the first $t$ label requests.
If $\PP(\operatorname{DIS}(V_t)) > 0$ for all values of $t$ in the algorithm, then
with probability $1$ the algorithm uses all $n$ label requests.
Technically, each claim below should be followed by the phrase,
``unless $\PP(\operatorname{DIS}(V_t))=0$ for some $t \leq n$,
in which case $\er(\hat{h}_n)=0$ so the bound trivially holds.''
However, to simplify the presentation, we will make this
special case implicit, and will not mention it further.

The high-level outline of this proof is to use $\PP(\operatorname{DIS}(V_t))$ as an
upper bound on $\sup_{h \in V_t} \er(h)$,
and then show $\PP(\operatorname{DIS}(V_t))$ is halved roughly every $\lambda=
\tilde{O}(\theta_f d)$ label requests. Thus, after
roughly $\tilde{O}(\theta_f d \log(1/\varepsilon))$ label requests,
any $h \in V_t$ should have $\er(h) \leq\varepsilon$.

Specifically, let $\lambda_n = \lceil8\theta_f (d \ln(8 e \theta
_f)+\ln(2n/\delta)) \rceil$.
If $n \leq\lambda_n$, the bound in the theorem statement trivially
holds, since the right-hand side exceeds $1$;
otherwise, consider some nonnegative $t \leq n-\lambda_n$ and
$t^{\prime} = t + \lambda_n$.
Let $X_{m_t}$ denote the point corresponding to the $t$th label request,
and let $X_{m_{t^\prime}}$ denote the point corresponding to label
request number $t^\prime$.
It must be that
\[
|\{X_{m_t+1},X_{m_t+2},\ldots,X_{m_{t^\prime}}\} \cap \operatorname{DIS}(V_t)| \geq
\lambda_n,
\]
which means there is an i.i.d. sample of size $\lambda_n$, with
distribution equivalent to
the conditional of $X$ given $\{X\in \operatorname{DIS}(V_t)\}$, contained
in $\{X_{m_t +1},\ldots,X_{m_{t^\prime}}\}$: namely, the first
$\lambda_n$ points in this subsequence
that are in $\operatorname{DIS}(V_t)$.

Now recall that, by classic results from the passive learning
literature (e.g., \cite{anthony99}),
this implies that on an event $E_{\delta,t}$ holding with probability
$1-\delta/n$,
\[
\sup_{h \in V_{t^\prime}} \er(h|\operatorname{DIS}(V_t)) \leq2\frac{d \ln({2
e \lambda_n}/{d}) + \ln({2 n}/{\delta})}{\lambda_n}.
\]
Also note that $\lambda_n$ was defined (with express purpose) so that
\[
2\frac{d \ln({2 e \lambda_n}/{d}) + \ln({2 n}/{\delta
})}{\lambda_n} \leq1/(2\theta_f).
\]
Recall that, since $\er(f) = 0$, we have $\er(h) = \PP(h(X) \neq f(X))$.
Since $f \in V_{t^\prime} \subseteq V_t$, this means for any $h \in
V_{t^{\prime}}$ we have $\{x \dvtx h(x) \neq f(x)\} \subseteq \operatorname{DIS}(V_t)$,
and thus
\begin{eqnarray*}
\sup_{h \in V_{t^{\prime}}} \PP\bigl(h(X) \neq f(X)\bigr) &=& \sup_{h \in
V_{t^{\prime}}} \PP\bigl(h(X)\neq f(X) | X \in \operatorname{DIS}(V_t)\bigr) \PP(\operatorname{DIS}(V_t)) \\
&=& \sup_{h \in V_{t^{\prime}}} \er(h|\operatorname{DIS}(V_t)) \PP(\operatorname{DIS}(V_t)) \leq\PP
(\operatorname{DIS}(V_t)) / (2 \theta_f).
\end{eqnarray*}
So $V_{t^{\prime}} \subseteq B(f,\PP(\operatorname{DIS}(V_t)) / (2\theta_f))$, and
therefore by monotonicity of $\PP(\operatorname{DIS}(\cdot))$ and the definition of
$\theta_f$
\[
\PP(\operatorname{DIS}(V_{t^\prime})) \leq\PP\bigl(\operatorname{DIS}\bigl(B\bigl(f,\PP
(\operatorname{DIS}(V_t))/(2\theta_f)\bigr)\bigr)\bigr) \leq\PP(\operatorname{DIS}(V_t))/2.
\]
By a union bound, $E_{\delta,t}$ holds for every $t \in\{i \lambda_n
\dvtx i \in\{0,1,\ldots,\lfloor n / \lambda_n\rfloor- 1\}\}$
with probability $\geq1-\delta$.
On these events, if $n \geq\lambda_n \lceil\log_2 (1/\varepsilon)
\rceil$, then (by induction)
\[
\sup_{h \in V_n} \er(h) \leq\PP(\operatorname{DIS}(V_n)) \leq\varepsilon.
\]
Solving for $\varepsilon$ in terms of $n$ gives the result (with a slight
increase in constants due to relaxing the ceiling functions).
\end{pf}

%%% TODO: IF POSSIBLE, INCLUDE A LOWER BOUND HERE TOO (IN TERMS OF
%THETA)

%s4.2 ###
\subsection{Known results on convergence rates for agnostic active learning}

We will now describe the known results for agnostic active learning
algorithms, starting with Algorithm 1. The key to the potential convergence
rate improvements of Algorithm 1 is that, as the region of disagreement $R$
decreases in measure, the error difference $\er(h|R)-\er(h^\prime|R)$ of
any classifiers $h,h^\prime\in V$ under the \textit{conditional}
sampling distribution (given $R$) can become significantly larger [by
a factor of $\PP(R)^{-1}$] than $\er(h)-\er(h^\prime)$, making it
significantly easier to determine which of the two is worse using a
sample of labeled data. In particular, \cite{hanneke07b} developed
a technique for analyzing this type of algorithm, and adapting that
analysis to the above definition of Algorithm 1 results in the following
guarantee.
\begin{theorem}[\cite{hanneke07b}]
Let $\hat{h}_n$ be the classifier returned by Algorithm 1 when allowed $n$
label requests, using the bounds (\ref{eqn:bbl-bounds}) and confidence
parameter $\delta\in(0,1/2)$. Then there exists a finite universal
constant $c$ such that, with probability $\geq1-\delta$, $\forall n
\in\nats$,
\begin{eqnarray*}
\er(\hat{h}_n) - \nu
&\leq& c \sqrt{\frac{\nu^2 \theta^2 (d \log
n + \log({1}/{\delta}))\log(({n+2\nu\theta})/({\nu\theta
}))}{n}} \\
&&{} + 2 \exp\biggl\{ -\frac{n}{c \theta^2 (d \log\theta+
\log({n}/{\delta}))} \biggr\}.
\end{eqnarray*}
\end{theorem}

Similarly, the key to improvements from Algorithm 2 is that as the
number $m$
of processed unlabeled data points increases, we only need to request
the labels of those data points in the region of disagreement of the
set of classifiers with near-optimal empirical error rates. Thus, if
the region of disagreement of classifiers with excess error $\leq
\varepsilon$ shrinks as $\varepsilon$ shrinks, we expect the frequency of
label requests to shrink as $m$ increases. Since we are careful not
to discard the best classifier, and the excess error rate of a
classifier can be bounded in terms of the $\Delta_m$ function, we end
up with a bound on the excess error which is converging in $m$, the
number of \textit{unlabeled} data points processed, even though we
request a number of labels growing slower than $m$. When this
situation occurs, we expect Algorithm 2 will provide an improved convergence
rate compared to passive learning.
Dasgupta, Hsu and Monteleoni \cite{dasgupta07} prove the following
convergence rate guarantee.
\begin{theorem}[\cite{dasgupta07}]
Let $\hat{h}_n$ be the classifier returned by Algorithm 2 when allowed $n$
label requests, using the threshold (\ref{eqn:dhm-bound}), and
confidence parameter $\delta\in(0,1/2)$. Then there exists a finite
universal constant $c$ such that, with probability $\geq1-\delta$,
$\forall n \in\nats$,
\begin{eqnarray*}
\er(\hat{h}_n) - \nu &\leq& c \sqrt{\frac{\nu^2 \theta(d \log
(({n+2\nu\theta})/({\nu\theta})) + \log({1}/{\delta}))}{n}}\\
&&{} +
c \biggl(d + \log\frac{1}{\delta}\biggr) \exp \Biggl\{ -\sqrt
{\frac{n}{c \theta(d + \log({1/\delta}))}}\Biggr\}.
\end{eqnarray*}
\end{theorem}

Note that, among other changes, this bound improves the dependence on
the disagreement coefficient $\theta$, compared to the bound for
Algorithm 1. In both cases, for certain ranges of $\theta$, $\nu$ and $n$,
these bounds can represent significant improvements in the excess error
guarantees, compared to the corresponding guarantees possible for
passive learning. However, in both cases, when $\nu> 0$ these bounds
have an \textit{asymptotic} dependence on $n$ of $\tilde{\Theta
}(n^{-1/2})$, which is no better than the convergence rates achievable
by passive learning (e.g., by empirical risk minimization). Thus, there
remains the question of whether either algorithm can achieve asymptotic
convergence rates strictly superior to passive learning for
distributions with nonzero noise rates. This is the topic we turn to next.

%s4.3 ###
\subsection{Active learning under Tsybakov's noise conditions}

It is known that for most nontrivial $\C$, for any $n$ and $\nu> 0$,
for every active learning algorithm there is some distribution with
noise rate $\nu$ for which we can guarantee excess error no better
than $\propto\nu n^{-1/2}$ \cite{kaariainen06}; that is, the
$n^{-1/2}$ asymptotic dependence on $n$ in the above bounds matches the
corresponding minimax rate, and thus cannot be improved as long as the
bounds depend on $\D_{XY}$ only via $\nu$ (and $\theta$). Therefore,
if we hope to discover situations in which these algorithms have
strictly superior asymptotic dependence on $n$, we will need to allow
the bounds to depend on a more detailed description of the noise
distribution than simply the noise rate $\nu$.

As previously mentioned, one way to describe a noise distribution
using a more detailed parametrization is to use Tsybakov's noise
conditions (Condition \ref{con:tsybakov}). In the context of passive
learning, this allows one to describe situations in which~the rate of
convergence is between $n^{-1}$ and $n^{-1/2}$, even when $\nu> 0$.
This raises the natural question of how these active learning
algorithms perform when the noise distribution satisfies this
condition with finite $\mu$ and $\kappa$ parameter values. In many
ways, it seems active learning is particularly well-suited to exploit
these more favorable noise conditions, since they imply that as we
eliminate suboptimal classifiers, the diameter of the remaining set
shrinks; thus, for finite $\theta$ values, the region of disagreement
should also be shrinking, allowing us to focus the samples in a
smaller region and accelerate the convergence.

Focusing on the special case of learning one-dimensional threshold
classifiers under a certain uniform marginal
distribution, Castro and Nowak \cite{castro08} studied conditions
related to
Condition \ref{con:tsybakov}. In particular, they studied a
threshold-learning algorithm that, unlike the algorithms described
here, takes $\kappa$ as \textit{input}, and found its convergence rate
to be $\propto(\frac{\log
n}{n})^{{\kappa}/({2\kappa-2})}$ when $\kappa> 1$, and
$\exp\{-c n\}$ for some ($\mu$-dependent) constant $c$, when $\kappa=
1$. Note that this improves over the $n^{-{\kappa}/({2\kappa-1})}$
rates achievable in passive learning \cite{castro08,tsybakov04}.
Subsequently, Balcan, Broder and Zhang \cite{balcan07} proved an
analogous positive result
for higher-dimensional linear separators, and Castro and Nowak \cite
{castro08} additionally
showed a related result for boundary fragment classes (see below); in
both cases,
the algorithm depends explicitly on the noise parameters.
Later, in a personal communication, Langford and Castro
claimed that in fact Algorithm 1 achieves this rate (up to log factors)
for the one-dimensional thresholds problem, leading to speculation that perhaps
these improvements are achievable in the general case as well (under
conditions on the disagreement coefficient).
Castro and Nowak \cite{castro08} also prove that a value $\propto
n^{-{\kappa}/({2\kappa-2})}$ (or $\exp\{-c^\prime n\}$, for some
$c^\prime$, when $\kappa=1$) is also a \textit{lower bound} on the
minimax rate for the threshold learning problem. In fact, a similar
proof to theirs can be used to show this same
lower bound holds for any nontrivial $\C$. For completeness,
a proof of this more general result is included
in Appendix D of the supplementary material \cite{hanneke10bsup}.

Other than the few specific results mentioned above,
it was not previously known whether Algorithm 1 or Algorithm 2, or
indeed \textit
{any} active learning algorithm,
generally achieves convergence rates that exhibit these types of improvements.

%s4.4 ###
\subsection{Adaptive rates in active learning: Algorithm 1}
\label{subsec:adaptive-agnostic}

The above observations open the question of whether these algorithms,
or variants thereof,
improve this asymptotic dependence on $n$.
It turns out this is indeed possible. Specifically, we have the following
result for Algorithm 1.
\begin{theorem}
\label{thm:BBL-adaptive}
Let $\hat{h}_n$ be the classifier returned by Algorithm 1 when allowed $n$
label requests, using the bounds (\ref{eqn:bbl-bounds}) and confidence
parameter $\delta\in(0,1/2)$. Suppose further that $\D_{XY}$
satisfies Condition \ref{con:tsybakov}. Then there exists a finite
($\kappa$- and $\mu$-dependent)
constant $c$ such that, for any $n\in\nats$, with probability $\geq
1-\delta$,
\[
\er(\hat{h}_n) - \nu\leq\cases{
2 \cdot \exp\biggl\{- \dfrac{n}{c \theta^2 (d \log n + \log(1 / \delta
))}\biggr\}, &\quad when
$\kappa=1$,\vspace*{2pt}\cr
c \biggl(\dfrac{\theta^2 (d\log n + \log(1 / \delta)) \log n
}{n}\biggr)^{{\kappa}/({2\kappa-2})}, &\quad when
$\kappa>1$.}
\]
%
%%% TODO: define r_0
\end{theorem}
\begin{pf}%[Proof of Theorem \ref{thm:BBL-adaptive}]
%The case of $\theta=0$ clearly holds, so we will focus on the
%nontrivial case of $\theta\geq1$.
We will proceed by bounding the \textit{label complexity},
or size of the label budget $n$ that is sufficient to guarantee,
with high probability, that the excess error of the returned classifier
will be at most $\varepsilon$ (for arbitrary $\varepsilon> 0$); with this
in hand, we can simply bound the
inverse of the function to get the result in terms of a bound on excess error.

Throughout this proof (and proofs of later results in this paper), we
will make frequent
use of basic facts about $\er(h|R)$. In particular, for any classifiers
$h, h^{\prime}$ and
set $R \subseteq\X$, we have $\er(h) = \er(h|R)\PP(R) + \er(h|\X
\setminus R) \PP(\X\setminus R)$;
also, if $\{x \dvtx h(x) \neq h^{\prime}(x)\} \subseteq R$, we have
$\er(h|\X\setminus R) - \er(h^{\prime} | \X\setminus R) = 0$ and therefore
$\er(h) - \er(h^{\prime}) = (\er(h | R) - \er(h^{\prime} | R))\PP(R)$.

Note that,
by Lemma \ref{lem:uniform} and a union bound,
on an event of probability $1-\delta$, (\ref{eqn:uniform})~holds with
$\delta^\prime= \delta/n$
for every set $Q$, relative to the conditional distribution given its
respective $R$ set,
for any value of $n$. For the remainder of this proof, we assume that this
$1-\delta$ probability event occurs.
In particular, this means that for every $h \in\C$ and every $Q$ set
in the algorithm,
$\LB (h,Q,\delta/n) \leq \er(h|R) \leq \UB (h,Q,\delta/n)$, for the set
$R$ that $Q$
is sampled under.

Our first task is to show that we never remove the ``good'' classifiers
from $V$.
We only remove a classifier $h$ from $V$ if
$h^{\prime} = \arg\min_{h^{\prime} \in V} \UB (h^{\prime},Q,\delta/n)$
has $\LB (h,Q,\delta/n) > \UB (h^{\prime},Q,\delta/n)$.
Each $h \in V$ has $\{x \dvtx h(x) \neq h^{\prime}(x)\} \subseteq \operatorname{DIS}(V)
\subseteq R$,
so that
\[
\UB (h^{\prime},Q,\delta/n) - \LB (h,Q,\delta/n) \geq \er(h^{\prime}|R)
- \er(h|R) = \frac{\er(h^{\prime})-\er(h)}{\PP(R)}.
\]
Thus, for any $h \in V$ with $\er(h) \leq \er(h^{\prime})$,
$\UB (h^{\prime},Q,\delta/n) - \LB (h,Q,\delta/n) \geq \er(h^{\prime} |
R) - \er(h|R) = (\er(h^{\prime}) - \er(h)) / \PP(R) \geq0$,
so that on any given round of the algorithm, the set $\{ h \in V \dvtx
\er(h) \leq \er(h^{\prime})\}$ is not removed from $V$.
In particular, since we always have $\er(h^{\prime}) \geq\nu$,
by induction this implies the invariant $\inf_{h \in V} \er(h) = \nu$,
and therefore also
\begin{eqnarray*}
\forall t\qquad \er(h_t) - \nu &=& \er(h_t) - \inf_{h \in V}\er(h) \\
&=&
\Bigl(\er(h_t|R)-\inf_{h\in V}\er(h|R)\Bigr)\PP(R) \leq\beta_t,
\end{eqnarray*}
where again the second equality is due to the fact that $\forall h \in V$,
$\{x \dvtx h_t(x) \neq h(x)\}$ $\subseteq \operatorname{DIS}(V) \subseteq R$.
%and therefore $\er(h_t | \X\setminus R) - \er(h | \X\setminus R) = 0$.
We will spend the remainder of the proof bounding the size of $n$
sufficient to guarantee some $\beta_t \leq\varepsilon$.
In particular, similar to the proof of Theorem~\ref{thm:cal-upper}, we
will see that
as long as $\beta_t > \varepsilon$, we will halve $\PP(\operatorname{DIS}(V))$ roughly every
$\tilde{O}(\theta^2 d \varepsilon^{{2}/{\kappa}-2})$
label requests,
so that the total number of label requests before some $\beta_t \leq
\varepsilon$
is at most roughly $\tilde{O}(\theta^2 d \varepsilon^{
{2}/{\kappa}-2} \log(1/\varepsilon))$.

%then we are guaranteed to have some $t$ with $\er(h_t) - \nu\leq
%Thus, we can simply find a sufficiently large $n$ so that
%the algorithm must either reach Step $3$ at least $\log(1/\varepsilon)$
%times before halting,
%or else it reaches Step $3$ fewer than $\log(1/\varepsilon)$
%times but has some $t$ with $\beta_t \leq\varepsilon$ anyway.}% end
%ignore

Recalling the definition of $h^{[k]}$ (from Definition \ref
{def:global-disagreement-coefficient}), let
%
%e4 ###
%
\begin{equation}
\label{eqn:Vtheta}
V^{(\theta)} = \biggl\{h \in V \dvtx\limsup_{k\rightarrow\infty
} \PP\bigl(h(X) \neq h^{[k]}(X)\bigr) > \frac{\PP(R)}{2\theta}\biggr\}.
\end{equation}
Note that after step 7, if $V^{(\theta)} = \varnothing$, then
\begin{eqnarray*}
\PP(\operatorname{DIS}(V))&\leq&\PP\Bigl(\operatorname{DIS}\Bigl(\Bigl\{h \in\C\dvtx\limsup
_{k\rightarrow\infty} \PP\bigl(h(X) \neq h^{[k]}(X)\bigr) \leq
\PP(R)/(2\theta)\Bigr\}\Bigr)\Bigr) \\
%&= \lim_{k^\prime\rightarrow\infty} \PP(\operatorname{DIS}(
%) \leq\PP(R)/(2\theta)\}))\\
&=& \lim_{k^\prime\rightarrow\infty} \PP\biggl(\operatorname{DIS}
\biggl(\bigcap_{k > k^\prime} B\bigl(h^{[k]},\PP(R)/(2\theta)
\bigr)\biggr)\biggr)\\
&\leq&\lim_{k^\prime\rightarrow\infty} \PP\biggl(\bigcap
_{k > k^\prime} \operatorname{DIS}\bigl(B\bigl(h^{[k]}, \PP(R)/(2\theta)
\bigr)\bigr)\biggr)\\
%&\leq\liminf_{k\rightarrow\infty} \PP(\operatorname{DIS}(\{h
&\leq&\liminf_{k \rightarrow\infty} \PP\bigl(\operatorname{DIS}
\bigl(B\bigl(h^{[k]}, \PP(R)/(2\theta)\bigr)\bigr)\bigr) \\
&\leq& \liminf
_{k\rightarrow\infty} \theta_{h^{[k]}} \frac{\PP(R)}{2\theta} =
\frac{\PP(R)}{2},
\end{eqnarray*}
so that we will %\ignore{definitely }
satisfy the condition in step 2 on the next round.
Here we have used the definition of $\theta$ in the final inequality
and equality.
On the other hand, if after step 7, we have $V^{(\theta)} \neq
\varnothing$, then
\begin{eqnarray*}
\varnothing&\neq&\biggl\{h \in V \dvtx\limsup_{k \rightarrow
\infty} \PP\bigl(h(X) \neq h^{[k]}(X)\bigr) > \frac{\PP(R)}{2\theta}\biggr\}
\\
& = &\biggl\{h \in V \dvtx\biggl(\frac{\limsup_{k \rightarrow
\infty} \PP(h(X) \neq h^{[k]}(X))}{\mu}\biggr)^{\kappa} >
\biggl(\frac{\PP(R)}{2\mu\theta}\biggr)^{\kappa}\biggr\}\\
& \subseteq &\biggl\{h \in V \dvtx\biggl(\frac{\diam(\er(h)-\nu;\C)}{\mu
}\biggr)^{\kappa} > \biggl(\frac{\PP(R)}{2\mu\theta}
\biggr)^{\kappa}\biggr\}\\
& \subseteq &\biggl\{h \in V \dvtx \er(h)-\nu> \biggl(\frac{\PP(R)}{2\mu
\theta}\biggr)^{\kappa}\biggr\}\\
& = &\Bigl\{h \in V \dvtx \er(h|R)-\inf_{h^\prime\in V} \er(h^\prime| R) >
\PP(R)^{\kappa-1}(2\mu\theta)^{-\kappa}\Bigr\}\\
& \subseteq &\Bigl\{h \in V \dvtx \UB (h,Q,\delta/n)-\min_{h^\prime\in V}
\LB (h^\prime,Q,\delta/n) > \PP(R)^{\kappa-1}(2\mu\theta)^{-\kappa
}\Bigr\}\\
& \subseteq &\Bigl\{h \in V \dvtx \LB (h,Q,\delta/n)-\min_{h^\prime\in V}
\UB (h^\prime,Q,\delta/n) \\
&&\hspace*{40.2pt} > \PP(R)^{\kappa-1}(2\mu\theta)^{-\kappa} -
4G(|Q|,\delta/n)\Bigr\}.
\end{eqnarray*}
Here, the third line follows from the fact that $\er(h^{[k]}) \leq
\er(h)$ for all sufficiently large $k$,
the fourth line follows from Condition \ref{con:tsybakov}, and the
final line follows from the definition of $\UB $ and $\LB $.
By definition, every $h \in V$ has $\LB (h,Q,\delta/n) \leq\min
_{h^\prime\in V} \UB (h^\prime,Q,\delta/n)$, so
for this last set to be nonempty after step 7, we must have $\PP
(R)^{\kappa-1}(2\mu\theta)^{-\kappa} < 4 G(|Q|,\delta/n)$.

%Thus, on any round for which we do not reach Step $3$, there is some
%$h \in V$ with
%$\PP(R)^{\kappa-1} (2\mu\theta)^{-\kappa} < \er(h | R) - \inf_{h^\prime
Combining these two cases ($V^{(\theta)} = \varnothing$ and $V^{(\theta
)} \neq\varnothing$),
since $|Q|$ gets reset to $0$ upon reaching step 3, we have that after
every execution of step 7,
%
%e5 ###
%
\begin{equation}
\label{eqn:PG-bound-1}
\PP(R)^{\kappa-1}(2\mu\theta)^{-\kappa} < 4 G(|Q|-1,\delta/n).
\end{equation}

If\vspace*{1pt} $\PP(R) \leq\frac{\varepsilon}{2 G(|Q|-1,\delta/n)} \leq\frac
{\varepsilon}{2 G(|Q|,\delta/n)}$,
then certainly $\beta_t \leq\varepsilon$ (by definition of $\beta_t
\leq2 G(|Q|,\delta/n) \PP(R)$).
So on any round for which $\beta_t > \varepsilon$, we must have
%
%e6 ###
%
\begin{equation}
\label{eqn:PG-bound-2}
\frac{\varepsilon}{2 G(|Q|-1,\delta/n)} < \PP(R).
\end{equation}
Combining (\ref{eqn:PG-bound-1}) and (\ref{eqn:PG-bound-2}), on any
round for which
$\beta_t > \varepsilon$,
%
%e7 ###
%
\begin{equation}
\label{eqn:PG-bound-3}
\biggl(\frac{\varepsilon}{2 G(|Q|-1,\delta/n)}\biggr)^{\kappa-1}(2\mu
\theta)^{-\kappa} < 4 G(|Q|-1,\delta/n).
\end{equation}
Solving for $G(|Q|-1,\delta/n)$ reveals that when $\beta_t > \varepsilon$,
%
%e8 ###
%
\begin{equation}
\label{eqn:Geps-bound}
4^{-1/\kappa} \biggl(\frac{\varepsilon}{2}\biggr)^{({\kappa
-1})/{\kappa}}(2\mu\theta)^{-1} < G(|Q|-1,\delta/n).
\end{equation}
Basic algebra shows that when $n \geq|Q| > d$, we have
\[
G(|Q|-1,\delta/n) \leq3 \sqrt{\frac{\ln({4/\delta}) + (d+1)
\ln(n)}{|Q|}}.
\]
Combining this with (\ref{eqn:Geps-bound}), solving for $|Q|$ and
adding $d$ to handle the case $|Q| \leq d$, we have that on any round
for which $\beta_t > \varepsilon$,
%
%e9 ###
%
\begin{equation}
\label{eqn:Q-bound}
|Q| \leq\biggl(\frac{2}{\varepsilon}\biggr)^{({2\kappa-2})/{\kappa
}}(6\mu\theta)^2 4^{2/\kappa} \biggl(\ln\frac{4}{\delta}+(d+1)\ln
(n)\biggr)+d.
\end{equation}
Since $\beta_t \leq\PP(R)$ by definition, and $\PP(R)$ is at least
halved each time we reach step~3,
we need to reach step 3 at most $\lceil\log_{2}(1/\varepsilon) \rceil$
times before we are guaranteed some $\beta_t \leq\varepsilon$. Thus, any
%
%e10 ###
%
\begin{equation}
\label{eqn:label-complexity}\quad
n \geq1 + \biggl(\biggl(\frac{2}{\varepsilon}\biggr)^{({2\kappa
-2})/{\kappa}}(6\mu\theta)^2 4^{2/\kappa} \biggl(\ln\frac{4}{\delta
}+(d+1)\ln(n)\biggr)+d\biggr)\log_2\frac{2}{\varepsilon}
\end{equation}
suffices to guarantee either some $|Q|$ exceeds (\ref{eqn:Q-bound}) or
we reach step 3 at least $\lceil\log_2 (1/\varepsilon)\rceil$ times,
either of which implies the existence of some $\beta_t \leq\varepsilon$.
The stated result now follows by basic inequalities
to bound the smallest value of $\varepsilon$ satisfying (\ref
{eqn:label-complexity})
for a given value of $n$.
\end{pf}

If the disagreement coefficient is finite, Theorem
\ref{thm:BBL-adaptive} can
often represent a significant improvement in convergence rate compared to
passive learning, where we typically expect rates of order $n^{-\kappa
/ (2\kappa-1)}$
\cite{mammen99,tsybakov04,castro08}; this gap is especially notable
when the disagreement
coefficient and $\kappa$ are small. Furthermore, the
bound matches (up to logarithmic factors)
the form of the minimax rate \textit{lower bound}
proved by Castro and Nowak \cite{castro08} for threshold classifiers
(where $\theta= 2$);
as mentioned, that lower bound proof can be generalized to any
nontrivial $\C$
(see Appendix D of the supplementary material \cite{hanneke10bsup}),
so that the rate of Theorem \ref{thm:BBL-adaptive}
is nearly minimax optimal for any nontrivial $\C$ with \textit{bounded}
disagreement coefficients. Also note that, unlike the upper bound analysis
of Castro and Nowak \cite{castro08}, we do not require the algorithm
to be given any extra
information about the noise distribution, so that this result is
somewhat stronger;
it is also more general, as this bound applies to an arbitrary
hypothesis class.

A refined analysis and minor tweaks to the algorithm should be able to
reduce the log factors in this result. For instance, defining $\UB$  and
$\LB$
using the uniform
convergence bounds of Alexander \cite{alexander84}, and using a
slightly more complicated algorithm closer to
the original definition \cite{balcan06,hanneke07b}---taking
multiple samples between bound evaluations,
allowing a larger confidence argument to the $\UB$  and $\LB $
evaluations---the $\log^2 n$ factor should reduce at least to $\log n \log\log n$,
if not further.
Also, as previously mentioned, it is possible to replace the quantities
$\PP(R)$ and $\PP(\operatorname{DIS}(V))$
in Algorithm 1 with estimators of these quantities based on a finite sample
of unlabeled data points,
while preserving the results of Theorem \ref{thm:BBL-adaptive} up to
constant factors.
We include an example of such estimators in Appendix C of the
supplementary material \cite{hanneke10bsup},
along with a sketch of how to modify the proof of Theorem
\ref{thm:BBL-adaptive}
to compensate for using these estimated probabilities.

%Furthermore, with high probability, using these estimators instead of
%exact probabilities
%only increase the total number of unlabeled data points used by the
%algorithm by a polynomial function.

%s4.5 ###
\subsection{Adaptive rates in active learning: Algorithm 2}
\label{subsec:dhm-adaptive-rates}

%In some sense, Theorem \ref{thm:BBL-adaptive} is somewhat surprising,
%since the bounds $\UB $ and $\LB $ used to define
%the set $V$ and the bounds $\beta_t$ are not themselves adaptive to
%the noise conditions.
%
Note that, as before, $n$~gets divided by $\theta^2$ in the rates
achieved by Algorithm 1.
As before, it is not clear whether any modification to the definitions
of $\UB $ and $\LB $
can reduce this exponent on $\theta$ from $2$ to $1$. As such, it is
natural to
investigate the rates achieved by Algorithm 2 under Condition
\ref{con:tsybakov}; we know that it
does improve the dependence on $\theta$ for the
worst case rates over distributions with any given noise rate, so we
might hope that it does the same
for the rates over distributions with any given values of $\mu$ and
$\kappa$. Unfortunately, we do
not presently know whether the original definition of Algorithm 2 achieves
this improvement. However, we
now present a slight modification of the algorithm, and prove that it
does indeed provide the desired
improvement in dependence on $\theta$, while maintaining the
improvements in the asymptotic dependence
on $n$. Specifically, consider the following definition for the
threshold in Algorithm 2:
%& \UB (h,Q,\delta^\prime) = \er_{Q}(h) + \hat{\bound}_{\C}(Q,\delta^
%& \LB (h,Q,\delta^\prime) = \er_{Q}(h) - \hat{\bound}_{\C}(Q,\delta^
%
%e11 ###
%
\begin{equation}
\label{eqn:dhm-tight}
\Delta_{m}\bigl(\LL,Q,h^{(y)},h^{(-y)},\delta\bigr) = 3 \hat{\bound}_{\C
}(\LL\cup Q,\delta;\LL),
\end{equation}
where $\hat{\bound}_{\C}(\cdot,\cdot;\cdot)$ is defined in
the \hyperref[app:bound]{Appendix}, based on a notion of local Rademacher
complexity studied by Koltchinskii \cite{koltchinskii06}. In
particular, the
quantity $\hat{\bound}_{\C}$ is known to be adaptive to Tsybakov's
noise conditions, in the sense that more favorable noise conditions
yield smaller values of $\hat{\bound}_{\C}$.
Using this definition, we have the
following theorem; due to space limitations, its proof is not presented
here, but is included in Appendix B of the supplementary material
\cite{hanneke10bsup}.
\begin{theorem}
\label{thm:tight-agnostic}
Suppose $\hat{h}_n$ is the classifier returned by Algorithm 2 with threshold
as in (\ref{eqn:dhm-tight}), when allowed $n$ label requests and given
confidence parameter $\delta\in(0,1/2)$. Suppose further that $\D_{XY}$
satisfies Condition \ref{con:tsybakov} with finite parameter values
$\kappa$ and $\mu$. Then there exists a finite ($\kappa$ and
$\mu$-dependent) constant $c$ such that, with probability $\geq1-\delta$,
$\forall n \in\nats$,
\[
\er(\hat{h}_n) - \nu\leq\cases{
\displaystyle c \biggl(d+\log\frac{1}{\delta}\biggr) \cdot \exp\Biggl\{- \sqrt
{\frac{n}{c \theta(d + \log(1/\delta))}}\Biggr\}, &\quad when
$\kappa=1$,\cr
\displaystyle c \biggl(\frac{\theta(d \log n + \log(1 / \delta))}{n}
\biggr)^{{\kappa}/({2\kappa-2})}, &\quad when
$\kappa>1$.}
\]
\end{theorem}

Note that this does indeed improve the dependence on $\theta$,
reducing its exponent from $2$ to $1$;
we do lose some in that there is now a square root in the exponent of
the $\kappa=1$ case;
however, as with Theorem \ref{thm:BBL-adaptive}, it is likely that
slight refinements to the definition
of $\Delta_m$ would reduce this (though we may also need to weaken the
theorem statement
to hold for any single $n$, rather than simultaneously for all $n$).

The bound in Theorem \ref{thm:tight-agnostic} is stated in terms of
the VC dimension $d$. However, for certain nonparametric hypothesis
classes, it is sometimes preferable to quantify the complexity of the
class in terms of a constraint on the \textit{entropy}
of the class, relative to the distribution $\D_{XY}$
(see e.g., \cite{vanderVaart96,tsybakov04,koltchinskii06,castro08}).
Specifically, for $\varepsilon\in[0,1]$, define
\[
\omega_{\C}(m,\varepsilon)= \E\mathop{\sup_{h_1,h_2 \in\C:}}_{\PP
\{h_1(X)\neq h_2(X)\}\leq\varepsilon}\bigl|\bigl(\er(h_1)-\er_{m}(h_1)\bigr)
-\bigl(\er(h_2)-\er_{m}(h_2)\bigr)\bigr|.
\]
\begin{condition}
\label{con:entropy}
There exist finite constants $\alpha> 0$ and $\rho\in(0,1)$ s.t.
$\forall m \in\nats$ and $\varepsilon\in[0,1]$,
$\omega_{\C}(m,\varepsilon) \leq\alpha\cdot\max\{\varepsilon
^{({1-\rho})/{2}} m^{-1/2}, m^{-{1}/({1+\rho})}\}$.
\end{condition}

In particular, the entropy with
bracketing condition used in the original minimax analysis
of Tsybakov \cite{tsybakov04} implies Condition \ref{con:entropy}
\cite{koltchinskii06},
as does the analogous condition for random entropy \cite
{gine03,gine06,koltchinskii08}. In
passive learning, it is known that empirical risk minimization
achieves a rate of order $n^{-\kappa/ (2\kappa+ \rho- 1)}$ under
Conditions~\ref{con:tsybakov} and \ref{con:entropy} \cite
{koltchinskii06,koltchinskii08}
(see also Appendix B of the supplementary material
\cite{hanneke10bsup}, especially (19) and Lemma 5),
and that this is sometimes minimax optimal \cite{tsybakov04}.
The following theorem gives a bound on the rate of convergence of the same
version of Algorithm~2 as in Theorem \ref{thm:tight-agnostic}, this
time in
terms of the entropy condition which, as before, is
faster than the passive learning rate when the disagreement
coefficient is finite. The proof of this result is included in
Appendix B of the supplementary material \cite{hanneke10bsup}.
\begin{theorem}
\label{thm:tight-agnostic-entropy}
Suppose $\hat{h}_n$ is the classifier returned by Algorithm 2 with threshold
as in (\ref{eqn:dhm-tight}), when allowed $n$ label requests and given
confidence parameter $\delta\in(0,1/2)$. Suppose further that $\D_{XY}$
satisfies Condition \ref{con:tsybakov} with finite parameter values
$\kappa$ and $\mu$, and Condition \ref{con:entropy} with parameter
values $\alpha$ and $\rho$. Then there exists a finite ($\kappa$,
$\mu$, $\alpha$ and $\rho$-dependent) constant $c$ such that, with
probability $\geq1-\delta$, $\forall n \in\nats$,
\[
\er(\hat{h}_n) - \nu\leq c \biggl(\frac{\theta\log(n / \delta
)}{n}\biggr)^{{\kappa}/({2\kappa+ \rho-2})}.
\]
\end{theorem}

Again, it is likely that refinements to the $\Delta_m$ definition may
lead to improvements in the log factor.
Also, although this result is stated for Algorithm 2, it is conceivable that,
by modifying Algorithm 1 to use
definitions of $V$ and $\beta_t$ based on $\hat{\bound}_{\C
}(Q,\delta;\varnothing)$, an analogous
result might be possible for Algorithm 1 as well.

It is worth mentioning that Castro and Nowak \cite{castro08} proved a
minimax lower bound for the
hypothesis class of \textit{boundary fragments}, with an exponent having
a similar
dependence on related definitions of $\kappa$ and $\rho$ parameters
to that of
Theorem \ref{thm:tight-agnostic-entropy}. Their result does provide a
valid lower bound
here; however, it is not clear whether their lower bound, Theorem
\ref{thm:tight-agnostic-entropy},
both, or neither is tight in the present context, since the value of
$\theta$ is not presently known for that particular problem, and
the matching upper bound of \cite{castro08} was proven under a
stronger restriction
on the noise than Condition \ref{con:tsybakov}.
However, see \cite{wang09} for an analysis of the disagreement
coefficient for other nonparametric hypothesis classes, characterized
by smoothness
of the decision surface.

% JOURNAL: include a paragraph discussing the following comment

%s5 ###
\section{Model selection}
\label{sec:aggregation}

While the previous sections address adaptation to the noise
distribution, they are still restrictive in that they deal with
hypothesis classes of limited expressiveness. That is, the assumption
of finite VC dimension implies a strong restriction on the variety of
classifiers one can
represent (or approximate) in the class; the entropy conditions allow slightly
more flexibility, but under nontrivial distributions, even the
entropy conditions imply a significant restriction on the
expressiveness of the class.
Thus, for algorithms restricted to classifiers from such a
restricted hypothesis class, it is often unrealistic to expect
convergence to the Bayes error rate.
We address this issue in this section by developing a
general algorithm for learning with a sequence of nested hypothesis
classes of increasing complexity, similar to the setting of Structural
Risk Minimization in passive learning \cite{vapnik82}. The objective
is to adapt, not only to the noise conditions, but also to the complexity
of the optimal classifier.
The starting point for this discussion is the assumption of a structure on
$\C$, in the form of a sequence of nested hypothesis classes:
\[
\C_1 \subset\C_2 \subset\cdots.
\]
Each class has an associated noise rate $\nu_i = \inf_{h \in\C
_i}\er(h)$, and we define $\nu_\infty= \lim_{i \rightarrow\infty}
\nu_i$. We also let $\theta_i$ and $d_i$ be the disagreement
coefficient and VC dimension, respectively, for the set $\C_i$.
We are interested in an algorithm that guarantees convergence in
probability of the error rate to $\nu_\infty$. We are particularly
interested in situations where $\nu_\infty= \nu^*$, a condition
which is realistic in this setting since the sets $\C_i$ can be
defined so that it is always satisfied, even while maintaining each
$d_i < \infty$
(see, e.g., \cite{devroye96}). Additionally, if we are so lucky as to
have some $\nu_i = \nu^*$, then we would like the convergence rate
achieved by the algorithm to be not significantly worse than running
one of the above agnostic active learning algorithms with hypothesis
class $\C_i$ alone.
%
%noise conditions, depending on whether it should be
%structure-dependent or structure-independent.}
In this context, we can define a structure-dependent version of
Tsybakov's noise condition as follows.
\begin{condition}
\label{con:tsybakov-dependent}
%(Structure-dependent Tsybakov noise)
For some nonempty $I \subseteq\nats$, for each $i \in I$, there exist finite
constants $\mu_i > 0$ and $\kappa_i \geq1$,
such that $\forall\varepsilon>0, \diam(\varepsilon; \C_i) \leq\mu_i
\varepsilon^{{1}/{\kappa_i}}$.
\end{condition}

Note that we do not require every $\C_i$, $i \in\nats$, to have
finite $\mu_i$ and $\kappa_i$, only some nonempty set $I \subseteq
\nats$;
this is important, since we might not expect $\C_i$ to satisfy
Condition \ref{con:tsybakov} for small indices $i$, where the
expressiveness is
quite restricted.

In passive learning, there are several methods for this type of
model selection which are known to preserve the convergence rates of each
class $\C_i$ under Condition \ref{con:tsybakov-dependent}
(e.g., \cite{tsybakov04,koltchinskii06}). In
particular, Koltchinskii \cite{koltchinskii06} develops a
method that performs this type of model selection; it turns out we can
modify Koltchinskii's method to suit our present needs in the context
of active learning; this results in a general active learning
model selection method that preserves the types of improved rates
discussed in the previous section. This modification is presented
below, based on using Algorithm 2 as a subroutine. (It should also be
possible to define an analogous method that uses Algorithm 1 as a subroutine
instead.)\vspace*{4pt}

\begin{bigboxit}
\textbf{Algorithm 3}\\
Input: nested sequence of classes $\{\C_i\}$, label budget $n$,
confidence parameter $\delta$\\
Output: classifier $\hat{h}_n$\\
{\vskip-2mm}\line(1,0){346}\\
%{\vskip-2mm}\line(1,0){377}\\
0. For $i = \lfloor\sqrt{n/2} \rfloor, \lfloor\sqrt{n/2} \rfloor
-1, \lfloor\sqrt{n/2} \rfloor-2, \ldots, 1$\\
1. \quad Let $\LL_{in}$ and $Q_{in}$ be the sets returned by Algorithm 2 run
with $\C_i$ and the \\
\qquad threshold (\ref{eqn:dhm-tight}), allowing
$\lfloor n / (2i^2) \rfloor$ label requests, and confidence $\delta/
(2i^2)$\\
2. \quad Let $h_{in} \leftarrow\mbox{\textsc{Learn}}_{\C_i}(\bigcup_{j \geq
i} \LL_{jn}, Q_{in})$\\
3. \quad If $h_{in} \neq\varnothing$ and $\forall j$ s.t. $i < j \leq
\lfloor\sqrt{n/2}\rfloor$,
%{\hskip2.3cm}
\phantom{aaaaaaaa}$\er_{\LL_{jn} \cup Q_{jn}}(h_{in}) - \er_{\LL_{jn}
\cup Q_{jn}}(h_{jn}) \leq\frac{3}{2} \hat{\bound}_{\C_j}\bigl(\LL_{jn}
\cup Q_{jn}, \delta/ (2j^2); \LL_{jn}\bigr)$\\
4. \quad\phantom{ai} $\hat{h}_n \leftarrow h_{in}$\\
5. Return $\hat{h}_n$
\end{bigboxit}

\vspace*{6pt}

The function $\hat{\bound}_{\cdot}(\cdot,\cdot;\cdot)$ is defined
in the \hyperref[app:bound]{Appendix}. This method can be shown to have a
confidence bound on its error rate converging to $\nu_{\infty}$ at a
rate never significantly worse than the original passive learning
method of Koltchinskii
\cite{koltchinskii06}, as desired. Additionally, we have the
following guarantee on the rate of convergence under Condition
\ref{con:tsybakov-dependent}.
%
% TODO: if I decide to include a few sentences explaining the high
%level proof, I should
% move the following sentence to after the theorem statement
%
The proof is similar in style to Koltchinskii's original proof, though
some care is needed due to the altered sampling distribution and the
constraint set $\LL_{jn}$. The proof is included in Appendix B of the
supplementary material~\cite{hanneke10bsup}.
\begin{theorem}
\label{thm:structure-dependent}
Suppose $\hat{h}_n$ is the classifier returned by Algorithm 3, when
allowed $n$ label requests and confidence parameter $\delta\in
(0,1/2)$. Suppose further that $\D_{XY}$ satisfies Condition
\ref{con:tsybakov-dependent}. Then there exist finite ($\kappa_i$ and $\mu
_i$-dependent) constants $c_i$ such that, with probability $\geq
1-\delta$, $\forall n \in\nats$,
\begin{eqnarray*}
&&
\er(\hat{h}_n) - \nu_\infty \\
&&\qquad\leq 3 \min_{i \in I} (\nu_i - \nu
_\infty) \\
&&\qquad\quad{}+
\cases{\displaystyle
c_i \biggl(d_i + \log\frac{1}{\delta}\biggr) \cdot
\exp\Biggl\{ - \sqrt{\frac{n}{c_i \theta_i (d_i + \log(
{1}/{\delta}))}}\Biggr\} ,
&\quad if $\kappa_i = 1$,\vspace*{2pt}\cr
\displaystyle c_i \biggl(\frac{\theta_i (d_i \log n + \log({1}/{\delta
}))}{n}\biggr)^{{\kappa_i}/({2\kappa_i-2})} ,&\quad if $\kappa_i
> 1$.}
\end{eqnarray*}
\end{theorem}

%%% TODO: can I include a few sentences explaining the proof at a high
%level?
In particular, if we are so lucky as to have $\nu_i = \nu^*$ for some
finite $i$,
then the above algorithm achieves a convergence rate not significantly worse
than that guaranteed by Theorem \ref{thm:tight-agnostic} for applying
Algorithm 2 directly,
with hypothesis class $\C_i$. Note that the algorithm itself has no
dependence on
the set $I$, nor has it any dependence on each class's complexity
parameters $d_i, \kappa_i, \mu_i, \theta_i$;
the adaptive behavior of the data-dependent bound $\hat{\bound}_{\C_j}$
allows\vspace*{1pt} the algorithm to adaptively ignore the returned classifier from
the runs of Algorithm 2 for
which convergence is slow, thus automatically selecting an index for
which the error rate is
relatively small.

As in the previous section, we can also show a variant of this result
when the complexities are quantified in terms of the entropy.
Specifically, consider the following condition and
theorem; the proof is in Appendix B of the supplementary material
\cite{hanneke10bsup}. Again, this
represents an improvement over known results for passive learning when
the disagreement coefficients are finite.
\begin{condition}
\label{con:entropy-aggregation}
For each $i \in\nats$, there exist finite constants $\alpha_i > 0$,
$\rho_i \in(0,1)$ s.t. $\forall m \in\nats$ and $\varepsilon\in[0,1]$,
$\omega_{\C_i}(m,\varepsilon) \leq\alpha_i \cdot\max\{
\varepsilon^{({1-\rho_i})/{2}} m^{-1/2}, m^{-{1}/({1+\rho
_i})}\}$.
\end{condition}
\begin{theorem}
\label{thm:structure-dependent-entropy}
Suppose $\hat{h}_n$ is the classifier returned by Algorithm 3, when
allowed $n$ label requests and confidence parameter $\delta\in(0,1/2)$.
Suppose further that $\D_{XY}$ satisfies
Conditions \ref{con:tsybakov-dependent} and \ref{con:entropy-aggregation}.
Then there exist finite
($\kappa_i$, $\mu_i$, $\alpha_i$ and $\rho_i$-dependent) constants
$c_i$ such that, with probability $\geq1-\delta$, $\forall n \in
\nats$,
\[
\er(\hat{h}_n) - \nu_\infty\leq3 \min_{i \in I} (\nu_i - \nu
_\infty) +
c_i\biggl(\frac{\theta_i \log(n/\delta)}{n}\biggr)^{{\kappa
_i}/({2\kappa_i + \rho_i - 2})}.
\]
\end{theorem}

In addition to these theorems for this structure-dependent version of
Tsybakov's noise conditions,
we also have the following result for a structure-independent noise
condition, in the sense that the
noise condition does not depend on the particular choice of $\C_i$
sets, but only on the distribution $\D_{XY}$
(and in some sense, the full class $\C= \bigcup_i \C_i$);
it may be particularly useful when the class $\C$ is universal, in the
sense that it can approximate any classifier.
\begin{theorem}
\label{thm:structure-independent}
Suppose the sequence $\{\C_i\}$ is constructed so that $\nu_{\infty}
= \nu^*$,
and $\hat{h}_n$ is the classifier returned by Algorithm 3, when allowed
$n$ label requests and confidence parameter $\delta\in(0,1/2)$.
Suppose that
there exists a constant $\mu> 0$ s.t. for all measurable $h \dvtx\X
\rightarrow\Y$,
$\er(h) - \nu^* \geq\mu\PP\{h(X) \neq h^*(X)\}$. Then there exists a
finite ($\mu$-dependent)
constant $c$ such that, with probability $\geq1-\delta$, $\forall n
\in\nats$,
\[
\er(\hat{h}_n) - \nu^* \leq c \min_{i \in\nats} (\nu_i - \nu^*) +
\biggl(d_i + \log\frac{i}{\delta}\biggr) \cdot \exp \Biggl\{
-\sqrt{\frac{n}{c i^2 \theta_i (d_i + \log(i / \delta
))}}\Biggr\}.
\]
\end{theorem}

The condition $\nu_{\infty} = \nu^*$ is quite easy to satisfy:
for example, $\C_i$ could be axis-aligned decision trees of depth $i$,
or thresholded polynomials of degree $i$,
or multi-layer neural networks with $i$ internal units, etc.
As for the noise condition in Theorem \ref{thm:structure-independent},
this would be satisfied whenever $\PP(|\eta(X) - 1/2| \geq c) = 1$
for some constant $c \in(0,1/2]$.
The case where $\er(h) - \nu^* \geq\mu\PP\{h(X) \neq h^*(X)\}^\kappa
$ for $\kappa> 1$ can be studied analogously, though the rate
improvements over passive learning are more subtle.

%rate no worse than (eqref passive bound from Koltchinskii), the rate
%achieved by the passive learning algorithm of Koltchinskii. The
%complete proof of this theorem is lengthy and fairly technical, so
%instead we provide a brief explanation. The bounds are roughly the
%same as the bounds from the previous theorem, and so Lepsky's
%thresholding technique adds at most an amount proportional to the
%bound value anyway, so that's not much loss.}
%
%technique for model selection by performing a sequence of hypothesis
%tests... Using this general technique, combined with a hypothesis test
%based on local Rademacher complexities, Koltchinskii (cite) developed
%a passive learning model selection algorithm that achieves the minimax
%rates for model selection in passive learning. In particular, that
%algorithm was shown by Koltchinskii (cite) to achieve a rate of
%convergence R(n,Dxy) $\propto$ ... under the structure-dependent
%definition of Tsybakov's noise conditions; thus, when the Bayes
%optimal classifier is contained in some class $\C_i$, this method
%conserves the convergence rate achievable by running a passive
%learning algorithm (in particular, Empirical Risk Minimization) with
%hypothesis class $\C_i$.}

% JOURNAL: include a paragraph discussion of the following comment
%necessarily nonverifiable.}

%s6 ###
\section{Conclusions}
\label{sec:conclusions}

Under Tsybakov's noise conditions, active learning can offer improved
asymptotic convergence rates compared to passive learning when the
disagreement coefficient is finite. It is also possible to preserve
these improved convergence rates when learning with a nested structure
of hypothesis classes, using an algorithm that adapts to both the noise
conditions and the complexity of the optimal classifier.

%%% TODO: I'd like to at least mention this entropy with bracketing
%thing
%%% Maybe it should wait for the Annals version, since that's the
%community that cares about this kind of thing anyway, moreso than the
%NIPS community which probably isn't too familiar with the idea.
%classes which have finite entropy with bracketing $\rho$, leading to
%rates similar to those above except that instead of having VC
%dimension factors in the bounds, the exponent becomes $\frac{\kappa}{2
%where the best known exponent is $\frac{\kappa}{2 \kappa+ \rho- 1}$,
%which is known to be tight in some cases; thus, for small disagreement
%coefficient, this represents an improved convergence rate.}% end ignore

\begin{appendix}\label{app:bound}
%s7 ###
\section*{Appendix: Definition of $\hat{\bound}$ and related quantities}

We define the quantity $\hat{\bound}_{\C}$ following Koltchinskii's
analysis of excess risk in terms
of local Rademacher complexity \cite{koltchinskii06}. The general
idea is to construct a bound
on the excess risk achieved by a given algorithm, such as empirical
risk minimization, via an application
of Talagrand's inequality. Such a bound should be based on a measure of
the expressiveness of the
set of functions $\C$; however, to bound the excess risk achieved by a
particular algorithm given a number of data
points, we need only measure the expressiveness of the set of functions
the algorithm is likely to select from.
For reasonable algorithms, such as empirical risk minimization, this
means the set of functions with reasonably small
excess risk. Thus, we can bound the excess risk of the algorithm in
terms of
a measure of expressiveness of the set of functions with relatively
small risk, typically referred to as
a \textit{local} complexity measure.
This reasoning is somewhat circular, in that first we must decide how
small to expect
the excess risk of the returned function to be before we can calculate
the local complexity measure,
which itself is used to calculate a bound on the risk of the returned
function. Thus, we define the
bound on the excess risk as a kind of fixed point. Furthermore, we can
estimate these quantities
using data-dependent confidence bounds, so that the excess risk bound
can be calculated without direct access to the
distribution. For the data-dependent measure of the expressiveness of
the function class, we can use
a Rademacher process. A detailed motivation and derivation can be found
in \cite{koltchinskii06}.

For our purposes, we add an additional constraint, by requiring the
functions we calculate the
complexity of to agree with the labels of a labeled set $\LL$. This is
helpful for us, since
given a set $Q$ of labeled data with true labels, for any two functions
$h_1$ and $h_2$ that
agree on the labels of $\LL$, it is always true that
$\er_{\LL\cup
Q}(h_1) - \er_{\LL\cup Q}(h_2)$ equals
the difference of the true empirical error rates. As we prove in the
supplement, as long as the set $\LL$ is chosen carefully (i.e., as in
Algorithm 2),
the addition of this constraint is essentially inconsequential, so that
$\hat{\bound}_{\C}$
remains a valid excess risk bound. The detailed
definitions are stated as follows.

For any function $f \dvtx\X\rightarrow\mathbb{R}$, and $\xi_1,\xi
_2,\ldots$ a sequence of independent random
variables with distribution uniform in $\{-1,+1\}$, define the \textit
{Rademacher process} for
$f$ under a finite set of (index, label) pairs $S \subset\nats\times
\Y$ as
\[
R(f;S) = \frac{1}{|S|}\sum_{(i,y) \in S} \xi_i f(X_i).
\]
The $\xi_i$ should be thought
of as internal variables in the learning algorithm, rather than being
fundamental to the learning problem.

For any two finite sets $\LL\subset\nats\times\Y$ and $S \subset
\nats\times\Y$,
define
\begin{eqnarray*}
\C[\LL] &=& \{h \in\C\dvtx \er_{\LL}(h) = 0\},\\
\hat{\C}(\varepsilon; \LL, S) &=& \Bigl\{h \in\C[\LL] \dvtx \er_{S}(h) - \min
_{h^\prime\in\C[\LL]} \er_{S}(h^\prime) \leq\varepsilon\Bigr\},\\
\hat{D}_{\C}(\varepsilon; \LL, S) &=& \sup_{h_1,h_2 \in\hat{\C
}(\varepsilon; \LL, S)} \frac{1}{|S|}\sum_{(i,y) \in S} \mathbh
{1}[h_1(X_i) \neq h_2(X_i)]
\end{eqnarray*}
and
\[
\hat{\phi}_{\C}(\varepsilon; \LL, S) = \frac
{1}{2}\sup_{h_1,h_2 \in\hat{\C}(\varepsilon; \LL, S)} R(h_1 - h_2; S).
\]
%
%%% TODO: can I name these things more descriptively? (i.e., can one of
%them be called \UB , another \bar{\UB }, etc.?)
For $\delta,\varepsilon> 0$, $m \in\nats$, define
$s_m(\delta) = \ln\frac{20 m^2 \log_2 (3 m)}{\delta}$ and
$\mathbb{Z}_{\varepsilon} = \{j \in\mathbb{Z} \dvtx2^j \geq\varepsilon
\}$,
and for any set $S \subset\nats\times\Y$, define the set $S^{(m)} =
\{(i,y) \in S \dvtx i \leq m\}$.
We use the following definitions from Koltchinskii
\cite{koltchinskii06} with only minor modifications.
\begin{definition}
\label{def:bound}
For $\varepsilon\in[0,1]$, and finite sets $S, \LL\subset\nats\times
\Y$, define
\[
\hat{U}_{\C}(\varepsilon,\delta; \LL, S) = \hat{K}\Biggl(\hat
{\phi}_{\C}(\hat{c}\varepsilon; \LL, S) + \sqrt{\frac
{s_{|S|}(\delta) \hat{D}_{\C}(\hat{c}\varepsilon; \LL, S)}{|S|}}
+\frac{s_{|S|}(\delta)}{|S|}\Biggr)
\]
and
\[
\hat{\bound}_{\C}(S, \delta; \LL) = \inf\Bigl\{
\varepsilon> 0 \dvtx\forall j \in\mathbb{Z}_{\varepsilon}, \min_{m
\in\nats} \hat{U}_{\C}\bigl(2^{j},\delta; \LL^{(m)}, S^{(m)}\bigr) \leq
2^{j-4}\Bigr\},
\]
where, for our purposes, we can take $\hat{K}=752$ %\ignore{$
and $\hat{c}=3/2$, %\ignore{and $\tilde{c}=3$},
though there seems to be room for improvement in these constants. %%%
%TODO: optimize the constants to get practical values for them
%discretization is left as a parameter. Here,
%for simplicity, we (arbitrarily) fix it to 2. }% end ignore
For completeness, we also define $\hat{\bound}_{\C}(\varnothing
,\delta; \C,\LL) = \infty$ by convention.
\end{definition}

%even if the label $y$ for some $(X_i,y) \in\LL$ disagrees with the
%true label $Y_i$, we can still calculate $\hat{\bound}_n(\delta;\C,
%and in particular we do not even need to know the true labels $Y_i$
%for the $X_i \in\LL_x$. This is because we only need
%to calculate \textit{differences} of empirical error rates, and every
%$h_1,h_2 \in\C[\LL]$ must \textit{agree}
%on the points in $\LL$ so that
%$\er_n(h_1) - \er_n(h_2) = \frac{1}{n} \sum_{i=1}^n \mathbh{1}[X_i
%Y_i])$.}% end ignore

%%% TODO: should I reformulate for the geometric definitions (inverse
%hats)? What's the bound improvement?

We will also define a related quantity, representing a
distribution-dependent version of $\hat{\bound}$, also
explored by Koltchinskii \cite{koltchinskii06}.
Specifically, for $\varepsilon> 0$, define
\[
\C(\varepsilon) = \{h \in\C\dvtx \er(h) - \nu\leq\varepsilon\}.
\]
For $m \in\nats$, let
\begin{eqnarray*}
\phi_{\C}(m,\varepsilon) &=& \E\sup_{h_1,h_2 \in\C(\varepsilon)}
\bigl|\bigl(\er(h_1) - \er_m(h_1)\bigr) - \bigl(\er(h_2) - \er_m(h_2)\bigr)\bigr|,\\
\tilde{U}_{\C}(m,\varepsilon,\delta) &=& \tilde{K}\Biggl(\phi_{\C
}(m,\tilde{c}\varepsilon) + \sqrt{\frac{s_m(\delta)\diam(\tilde
{c}\varepsilon;\C)}{m}} + \frac{s_m(\delta)}{m}\Biggr)
\end{eqnarray*}
and
\[
\tilde{\bound}_{\C}(m,\delta) = \inf\{\varepsilon
> 0 \dvtx\forall j \in\mathbb{Z}_{\varepsilon}, \tilde{U}_{\C
}(m,2^j,\delta) \leq2^{j-4}\},
\]
where, for our purposes, we can take $\tilde{K} = 8272$ and $\tilde
{c} = 3$.
For completeness, we also define $\tilde{\bound}_{\C}(0,\delta) =
\infty$.

%s7.1 ###
\subsection{Definition of $r_0$}
\label{subsec:r0}

In Definition \ref{def:disagreement-coefficient}, we took $r_0 = 0$.
If $\theta< \infty$, then this choice is usually relatively harmless.
However, in some cases, setting $r_0 = 0$ results in a suboptimal, or
even infinite, value of $\theta$, which is undesirable. In these
cases, we would like to set $r_0$ as large as possible while
maintaining the validity of the bounds. If we do this carefully enough, we
should be able to establish bounds that, even in the worst case when
$\theta= 1/r_0$, are never worse than the bounds for some analogous
passive learning method; however, to do this requires $r_0$ to depend
on the parameters of the learning problem: namely, $n$, $\delta$,
$\C$ and $\D_{XY}$.
The effect of a larger $r_0$ can sometimes be dramatic, as there are
scenarios where
$1 \ll\theta\ll1/r_0$~\cite{hanneke08a};
we certainly wish to distinguish between such scenarios, and those
where $\theta\propto1/r_0$.

Generally, depending on the bound we wish to prove, different values of
$r_0$ may be appropriate.
For the tightest bound in terms of $\theta$ proven in the Appendices
(namely, Lemma 7
of Appendix B in the supplementary material \cite{hanneke10bsup}),
the definition of $r_0 = r_{\C}(n,\delta)$ in (\ref{eqn:r0}) below
gives a good bound.
For the looser bounds (namely, Theorems \ref{thm:tight-agnostic}
and \ref{thm:tight-agnostic-entropy}),
a larger value of $r_0$ may provide better bounds;
however, this same general technique can be employed to define a good
value for $r_0$ in
these looser bounds as well, simply using upper bounds on (\ref
{eqn:r0}) analogous to how the theorems
themselves are derived from Lemma 7 in %\ignore{the proof of Lemma 8 in
%}
Appendix B \cite{hanneke10bsup}.
Likewise, one can state analogous refinements of $r_0$ for Theorems
\ref{thm:cal-upper}--\ref{thm:BBL-adaptive},
though for brevity these are left for the reader's independent consideration.
\begin{definition}
\label{def:r0}
Define
%
%e12 ###
%
\begin{equation}
\label{eqn:tildem}\qquad
\tilde{m}_{\C}(n,\delta) = \min\Biggl\{m \in\nats\dvtx n \leq
\log_2 \frac{4 m^2}{\delta} + 2e \sum_{\ell=0}^{m-1}\PP
(\operatorname{DIS}(\C(6\tilde{\bound}_{\C}(\ell,\delta)))
)\Biggr\}
\end{equation}
and
%
%e13 ###
%
\begin{equation}
\label{eqn:r0}
r_{\C}(n,\delta) = \max\Biggl\{\frac{1}{\tilde{m}_{\C}(n,\delta
)}\sum_{\ell=0}^{\tilde{m}_{\C}(n,\delta)-1} \diam(6\tilde
{\bound}_{\C}(\ell,\delta);\C), 2^{-n}\Biggr\}.
\end{equation}
\end{definition}

We use this definition of $r_0 = r_{\C}(n,\delta)$ in all of the main proofs.
In particular, with this definition, Lemma 7 of Appendix B
\cite{hanneke10bsup}
is never significantly worse than the analogous known result for
passive learning
(though it can be significantly better when $\theta\ll 1/r_0$).
%%% TODO: maybe I actually want to state the definition of r_0 for both
%of these bounds, which gives the fallbacks, so the reader doesn't have
%to.
\end{appendix}

\section*{Acknowledgments}

I extend my sincere thanks to Larry Wasserman for numerous helpful
discussions and also to John Langford for initially pointing out to me
the possibility of $A^2$ adapting to Tsybakov's noise conditions for
threshold classifiers. I would also like to thank the anonymous
referees for extremely helpful suggestions on earlier drafts.

\begin{supplement}[id=supp]
\stitle{Proofs and Supplements for ``Rates of Convergence in Active Learning''}
\slink[doi]{10.1214/10-AOS843SUPP}
\sdatatype{.pdf}
\sfilename{supplement-active-rates-annals.pdf}
\sdescription{The supplementary material contains three additional
Appendices, namely,
Appendices B, C and D.
Specifically, Appendix B provides detailed proofs of Theorems
%and \ref{thm:structure-independent},
\ref{thm:tight-agnostic}--\ref{thm:structure-independent},
as well as several abstract lemmas from which these results are derived.
Appendix C discusses the use of estimators in Algorithm 1.
Finally, Appendix D includes a proof of a general minimax lower bound
$\propto n^{-{\kappa}/({2\kappa-2})}$ for any nontrivial
hypothesis class, generalizing a result of Castro and Nowak \cite{castro08}.}
\end{supplement}

%suskaldyti doi

%
\printaddresses

\end{document}